\documentclass[11pt]{article}
\usepackage{amsmath,graphicx}
\usepackage{amsthm}
\usepackage{mathtools}
\usepackage{amsfonts}
\usepackage{hyperref}
\usepackage[title]{appendix}
\usepackage[font={small,sl},labelfont={sf,up}]{caption}

\allowdisplaybreaks

\newtheorem{theorem}{Theorem}[section]
\newtheorem{corollary}[theorem]{Corollary}
\newtheorem{lemma}{Lemma}[section]

\newtheorem{definition}{Definition}[section]
\newtheorem{claim}[theorem]{Claim}
\newtheorem{fact}[theorem]{Fact}

\newcommand*\N{\mathbb{N}}
\newcommand*\E{\mathbb{E}}
\newcommand*\R{\mathbb{R}}
\newcommand*\X{\mathrm{Del}_X}
\newcommand*\Y{\mathrm{Del}_Y}

\newcommand{\A}{\mathcal{A}}

\newcommand{\mcalE}{\mathcal{E}}
\newcommand{\F}{\mathcal{F}}

\newcommand{\Pa}{\mathcal{P}}

\newcommand{\bG}{\mathbb{G}}
\newcommand{\bP}{\mathbb{P}}
\newcommand{\bE}{\mathbb{E}}

\newcommand{\bF}{\mathbb{F}}
\newcommand{\bN}{\mathbb{N}}
\newcommand{\ep}{\varepsilon}

\title{The Normalized Matching Property in Random and Pseudorandom Bipartite Graphs}
\author{Niranjan Balachandran\footnote{Department of Mathematics, Indian Institute of Technology Bombay. Email: niranj@math.iitb.ac.in.} \ and Deepanshu Kush\footnote{Department of Computer Science, University of Toronto. Email: deepkush@cs.toronto.edu. Work was performed while the author was at IIT Bombay.}}

\begin{document}

\maketitle

\begin{abstract}
A simple generalization of the Hall's condition in bipartite graphs, the Normalized Matching Property (NMP) in a graph $G(X,Y,E)$ with vertex partition $(X,Y)$ states that for any subset $S\subseteq X$, we have $\frac{|N(S)|}{|Y|}\ge\frac{|S|}{|X|}$. In this paper, we show the following results about the Normalized Matching Property in random and pseudorandom graphs.
\begin{enumerate}

\item We establish $p=\frac{\log n}{k}$ as a sharp threshold for having NMP in $\mathbb{G}(k,n,p)$, which is the graph with $|X|=k,|Y|=n$ (assuming $k\le n\leq \exp(o(k))$), and in which each pair $(x,y)\in X\times Y$ is an edge  independently with probability $p$. This generalizes a classic result of Erd\H{o}s-R\'enyi on the $\frac{\log n}{n}$ threshold for having a perfect matching in  $\bG(n,n,p)$.
\item We also show that a pseudorandom bipartite graph - upon deletion of a vanishingly small fraction of vertices - admits NMP, provided it is not too sparse. More precisely, a bipartite graph $G(X,Y)$, with $k=|X|\le |Y|=n$,  is said to be Thomason pseudorandom (following A. Thomason (Discrete Math., 1989)) with parameters $(p,\ep)$ if each $x\in X$ has degree at least $pn$ and each pair of distinct $x, x'\in X$ has at most $(1+\ep)p^2n$ common neighbors. We show that for  any large enough $(p,\ep)$-Thomason pseudorandom graph $G(X,Y)$, there are ``tiny" subsets $\X\subset X, \  \Y\subset Y$ such that the subgraph $G(X\setminus \X,Y\setminus \Y)$ has NMP,  provided $p \gg\tfrac{1}{k}$. En route, we prove an ``almost'' vertex decomposition theorem: Every such Thomason pseudorandom graph admits - excluding a negligible portion of its vertex set - a partition of its vertex set into graphs that we call \emph{Euclidean trees}. These are trees that have NMP, and which arise organically through the Euclidean GCD algorithm. 
\end{enumerate}
\end{abstract}

\newpage

\section{Introduction}
Consider the following problems:
\begin{enumerate}
\item \label{prob:array} Suppose $k\le n$ are positive integers. By a {\it $k\times n$ star-array} (or simply star-array), we mean a $k\times n$ array whose entries are symbols from the set $\{0,\star\}$.  Given a $k\times n$ star-array, when is it possible to replace some of the $\star$ entries of the array by non-negative integers such that in the resulting array all the row sums  equal $R$, and all the column sums equal $C$ for some integers $R,C>0$? 
\item \label{prob:sumcayley} Let $q$ be a sufficiently large prime power and suppose $X, Y\subset\bF_q$ with $|Y|=10|X|$, $|X|\ge q/100$. Is it possible to label each element of $Y$ with some element of $X$ such that each element of $X$ appears as a label exactly $10$ times, and further, for each $y\in Y$ labeled $x$, the sum $x+y$ is a quadratic residue? More generally, one can ask the same question with a subgroup $H\subset\bF_q^*$ instead of the set of quadratic residues.
\end{enumerate}

In both  the problems posed above, there is a natural bipartite graph $G(X,Y,E)$ that captures the problem in its essence: Given a star-array $\A$, let $X$ and $Y$ denote the set of rows and columns of $\A$ respectively, and a vertex $x\in X$ is adjacent to $y\in Y$ in $G$ if and only if the $(x,y)$ entry of  $\A$ corresponding to a $\star$. For the second problem consider the bipartite graph $G(X,Y,E)$ where $X, Y$ are the given sets, and the pair $(x,y)$ is an edge in $G$ if and only if $x+y\in H$. 

 In the rest of the paper,  $G(X,Y)$  shall denote a bipartite graph with vertex partition $(X,Y)$; we shall drop the $E$ in our notation for convenience. We say that $G=G(X,Y)$ has the Normalized Matching Property (NMP for short) if: For any $S\subseteq X$, if we denote by $N(S)$, its set of neighbors in $Y$, then $\frac{|N(S)|}{|Y|}\ge\frac{|S|}{|X|}$. In particular, if $|X|=|Y|$, then this is the familiar Hall's condition for the existence of a perfect matching in $G$.  
 
 The following theorem of Kleitman  \cite{kleit} gives us an equivalent formulation of NMP in bipartite graphs:
 \begin{theorem} The following statements are equivalent:
 \begin{itemize} \label{thm:kleitman}
 \item $G$ with $|X|=k,|Y|=n$ has NMP.
 \item  For any independent set $I$ in $G$, $\frac{|I_X|}{k}+\frac{|I_Y|}{n}\le 1$.
\item There exists a multiplicity function $m:E\to\N_0=\N\cup\{0\}$ such that $\displaystyle\sum_{\substack{e\ni x\\e\in E}} m(e)$ (resp. $\displaystyle\sum_{\substack{e\ni y\\e\in E}} m(e)$) is equal for all $x\in X$ (resp. for all $y\in Y$).\end{itemize}\end{theorem}

It is easy to see that the problems posed above simply ask if the associated bipartite graphs have NMP by virtue of the third part of Theorem \ref{thm:kleitman}.

%One may frame the desired property in terms of  this graph. For instance, if  the corresponding graph is biregular, then this property trivially holds. On the other hand, if all the vertices of $X$ have the same degree, then it is not necessarily the case that the desired property holds; indeed, it could happen that some vertex of $Y$ has no edge incident to it at all,  as the edges between $X$ and $Y$ may form a cluster with some proper subset of $Y$, say.
The Normalized Matching Property in bipartite graphs was introduced by Graham and Harper \cite{GH} and subsequently has been a focus of study in bipartite graphs in several papers (for instance \cite{kleit, WHD}) and some monographs as well (for instance \cite{And2, Engel}). The notion also extends very naturally to finite ranked posets; for a ranked poset $P$, let $L_i$ denote the set of all elements of $P$ with rank $i$. Then we say that $P$ has NMP if for each $i$, the bipartite graph of poset covering  relations between $L_i$ and $L_{i+1}$ has NMP. NMP posets are objects of great interest specifically in related decomposition problems (see \cite{Griggs,GWL, HLS, Wang, WZ} for some decompositions results). As a concrete instance, the Griggs conjecture which states that any unimodal NMP poset admits a nested chain decomposition (see \cite{HLS} or \cite{LC} for more details on what the definitions are) is still open - even for posets of rank $3$ - despite several attacks on the problem.

As it turns out, many interesting finite ranked posets arising from finite geometric structures have NMP. Indeed, the Boolean poset, the poset of affine flats in a finite projective $n$-dimensional space and the poset of the subgroup lattice of abelian $p$-groups all have NMP (see \cite{thom2,Wang, WZ} respectively), i.e., in each of these posets, the associated bipartite graphs on the sets of elements of successive ranks within these posets have NMP.  As is the case with Hall's theorem for bipartite graphs, it is clear that graphs with ``high density'' are more likely to possess NMP. But in each of the instances listed above, the associated bipartite graphs are very sparse. This raises the following natural question: \textit{At what density does a typical bipartite graph have NMP}? 

To formulate the above question more precisely, we set up some asymptotic terminology and notation.  Given functions $f,g$, we write $f\gg g$ (resp. $f\ll g$) if $\displaystyle\lim_{n\to\infty}\frac{f(n)}{g(n)}\to \infty$ (resp. $\frac{f(n)}{g(n)}\to 0$). We also write $f=o(g)$ to denote that $f\ll g$. We write $f=O(g)$ (resp. $f=\Omega(g)$)  if there exists an absolute constant $C>0$ and $n_0$ such that for all $n\ge n_0, |f(n)|\le C|g(n)|$ (resp. if $|f(n)|\ge C|g(n)|$). If the constant $C$ involves a related parameter $\ep$, then we write $f=O_{\ep}(g)$ (resp. $f=\Omega_{\ep}(g)$) to indicate the dependence of the implicit constant on the parameter $\ep$. 

	To formalize the question posed above, we recall some standard terminology from the theory of random graphs. For a probability space $(\Omega, \bP)$ we say that an event $\mcalE_n$ that depends on a parameter $n$ {\it occurs  with high probability} (abbreviated as {\it whp}) if $\bP(\mcalE_n)\to 1$ as $n\to\infty$.  A graph property $\Pa$ is simply a collection of graphs, and a graph property is called monotone if whenever $G\in\Pa$ and $G\subset H$ then $H\in\Pa$ as well. The Erd\H{o}s-R\'enyi random graph model $\bG(n,p)$ introduced in \cite{ER1} is the random graph where the vertex set is the set $[n]:=\{1,\ldots,n\}$ and each pair $\{i,j\}$ is an edge with probability $p=p(n)$ independently. A monotone graph property $\Pa$ is said to have a threshold $p_0=p_0(n)$ if whenever $p\gg p_0$ then $\mathbb{G}(n,p)$ has property $\Pa$ {\it whp}, and if  $p\ll p_0$ then {\it whp} $\mathbb{G}(n,p)$ does not have property $\Pa$.  A property $\Pa$ is said to have a {\it sharp} threshold $p_0(n)$ if for $\ep>0$ and $p\ge(1+\ep)p_0$, $\mathbb{G}(n,p)$ has property $\Pa$ {\it whp} and for $p\le (1-\ep)p_0$,   $\mathbb{G}(n,p)$ does not have property $\Pa$ {\it whp}.
 
The seminal paper of Erd\H{o}s and R\'enyi \cite {ER1} established sharp thresholds for several very natural monotone graph properties. A theorem of Bollobas and Thomason \cite{BT} showed that every monotone graph property admits a threshold. However, not all graph properties admit sharp thresholds; for instance, the property ``$\bG(n,p)$ contains a cycle'' admits a threshold which is sharp on one side but not the other (see \cite{rg} for more on sharp thresholds). In fact, the problem of determining sharp thresholds (if the graph property admits one) is a very popular motif in the theory of random graphs.

For bipartite graphs, Erd\H{o}s and R\'enyi also introduced the random bipartite model $\bG(n,n,p)$ where the vertex set is partitioned into two sets $X,Y$ of size $n$ each, and each pair $\{x,y\}$ with $x\in X, y\in Y$ is in $\bG(n,n,p)$ independently with probability $p$. One of the first results in this model is the result that $\frac{\log n}{n}$ is a sharp threshold for the existence of a perfect matching in $\bG(n,n,p)$ \cite{ER2}. As observed earlier, if $k=n$, NMP is the same as Hall's condition for bipartite graphs, so it is natural to seek the threshold for NMP in a slightly more general model for bipartite random graphs, which is what the question previously posed seeks to do. 

Suppose $k\le n$ are positive integers, and let $0\le p\le 1$. Let $\mathbb{G}(k,n,p)$ denote the random bipartite graph with the vertex partition given by $(X,Y)$ with $|X|=k,|Y|=n$, and each pair $(x,y)\in X\times Y$ is an edge in $\mathbb{G}$ independently with probability $p$. {Here both $k$ and $n$ should be thought of as parameters growing to infinity with $n$ being a function of $k$ that always satisfies $n\geq k$.} Our first main result in this paper establishes a sharp threshold for NMP in the sense stated above:
%To make the sense of probability more precise, we shall use the Landau asymptotic notation. . %For a parameter $\alpha>0$, we write $f=O_{\alpha}(g)$ (resp. $f=\Omega_{\alpha}(g)$) if the corresponding constant $C=C_{\alpha}$ depends on the parameter $\alpha$.
%More precisely: Given $\delta>0$, do there exist $k_0(\delta),n_0(\delta)$ and a threshold function $p=p(n,k)$ such that \begin{eqnarray*}\textrm{If\ } &p\gg p(n,k),& \mathbb{G}(k,n,p)\textrm{\  has\  NMP\  with\  probability\  at\  least\ } 1-\delta\\ 
%\textrm{If\ }&p\ll p(n,k),& \mathbb{G}(k,n,p)\textrm{\  does\ {\it not}\ have\   NMP\  with\  probability\  at\  least\ } 1-\delta\end{eqnarray*} for $k\ge k_0, n\ge n_0$? 
\begin{theorem}\label{thresh_NMP}
Suppose $k\leq n(k) \le \exp(o(k))$, and let $0<\ep,\delta <1$. There exists $k_0 = k_0 (\ep,\delta)$ such that for $k\ge k_0(\ep,\delta)$
\begin{enumerate}
\item If $p\geq \frac{(1+\ep)\log n}{k}$ then $\bP[\mathbb{G}(k,n,p)\textrm{\  has\  NMP}]\ge 1-\delta$. 
\item If $p\leq \frac{(1-\ep)\log n}{k}$ then $\bP[\mathbb{G}(k,n,p)\textrm{\  has\  NMP}]\le \delta$.
\end{enumerate}
\end{theorem}
In other words, $\mathbb{G}(k,n,p)$ has a sharp threshold for NMP at $p=\frac{\log n}{k}$.

 Note that if $n> \exp(k)$ or equivalently, if $\log n>k$, then the expression for our threshold exceeds one. Also, for each fixed $p<1$, if $C>1+\log(\frac{1}{1-p})$ and $n\ge\exp(Ck)$, then a simple computation shows that the probability that $Y$ has at least one isolated vertex is bounded away from zero (this will be clear from the proof of Theorem \ref{thresh_NMP}; see Lemma \ref{lem:ub}). Hence, the range for $n$ in the statement of the theorem is essentially the widest possible one if one seeks a sharp threshold.
   
Let us now return to the problems at the beginning of this section.  To check if a given bipartite graph has NMP is  computationally simple: form a bigger new bipartite graph $G'(X',Y')$ with $|X'|=|Y'|=nk$ with $X'$ consisting of by  $n$ copies of $X$, $Y'$consisting of $k$ copies of $Y$, and $x'y'$ being an edge in $G'$ if and only if $xy$ was an edge in $G$. Then it is straightforward to see that  $G$ has NMP if and only if $G'$ admits a perfect matching. Hence either problem admits a computationally simple solution. But let us relax our requirement and seek an  answer only in an {\it approximate} sense: For the first problem, is it possible to replace each $\star$ entry with a non-negative integer such that {\it with the exception of a negligible proportion of the rows/columns}, the remaining rows and columns satisfy the aforementioned property? Or in the second problem, can we ignore a negligible proportion of elements from both sets $X, Y$, so that the desired property holds for the remaining elements? Since either of the originally posed problems is equivalent to asking if a given bipartite graph has NMP, this approximate version asks if a given bipartite graph ``almost'' has NMP in a certain sense that we shall formalize below.

The bipartite graph considered in the second problem (with the subsets of $\bF_q$) possesses certain regularity properties that are best described as ``random-like'' - as we shall soon see. Taking a cue from this, we impose the following reasonable hypotheses on bipartite graphs that we shall consider: If all the vertices of $X$ have ``almost''  the same degree, and suppose that no two vertices of $X$ have ``too many'' common neighbors in $Y$ (so that there isn't a clustering of edges between some subsets of $X$ and $Y$), is there an affirmative answer to the approximate  version for these problems? 

  To formulate this in more precise terms, we need the notion of a \textit{pseudorandom bipartite graph}. The notion of pseudorandomness was first introduced by Thomason in the 80s  \cite{thom1} and pseudorandomness in graphs is a  well-studied notion (see \cite{KS} for a definitive survey). One of the more popular and  well-understood models for pseudorandomness in graphs is the notion of an $(n,d,\lambda)$ graph (see \cite{AS15}).  An $(n,d,\lambda)$ graph is a graph on $n$ vertices which is $d$-regular and which satisfies the following property: If $d=\lambda_1\ge\lambda_2\ge\cdots\ge\lambda_n$ are the eigenvalues of $G$ then $|\lambda_i|\le\lambda$ for all $i>1$. 
  
 Pseudorandom graphs, as the name suggests, have some properties very reminiscent of random graphs, and the most well-known is the Expander-Mixing Lemma (see \cite{AS15}): Suppose $G$ is an $(n,d,\lambda)$ graph. If $U,W\subset V(G)$ then $|e(U,W)-\frac{d |U||W|}{n}|\le\lambda\sqrt{|U||W|}$, where $e(U,W)$ denotes the number of edges of the form $uw$ with $u\in U$ and $w\in W$. 
  
As mentioned earlier, Thomason introduced the notion of pseudorandomness which is a little more general, and in particular,  we shall - in this paper - confine our attention to the notion of pseudorandomness in bipartite graphs as proposed by Thomason in \cite{thom2}.   
  %Secondly, this notion of pseudorandomness is easy to verify (in algorithmic terms) for a given graph and a specified value of $\ep$. 
\begin{definition}\label{def:pseudo}Suppose $0<p<1$, and $0\le\ep<1$. 
A bipartite graph $G$ with vertex classes $X$ and $Y$ of sizes $k$ and $n$ respectively with $k\leq n $ is called Thomason pseudorandom with parameters $(p,\varepsilon)$ if every vertex in $X$ has degree at least $pn$,  and every pair of distinct vertices in $X$ have at most $p^2n(1 +\varepsilon)$ neighbors in common.
\end{definition}
At this juncture, a few remarks  are in order. Thomason's original definition in \cite{thom2} actually only considers bipartite graphs with $|X|=|Y|=n$.  Secondly, Thomason's definition in \cite{thom2} is more in line with the original notion of pseudorandomness in \cite{thom1}: A graph $G(X,Y)$ is  pseudorandom with parameters $(p,\mu)$ for some $\mu\geq 0$ where the second condition  states that every pair of vertices in $X$ have at most $p^2n+\mu$ common neighbors. The definition that we shall be using is a relaxation of the restriction that $|X|=|Y|$, but also a restriction to the more natural and intuitive case where $\mu\le \ep p^2n$. 
 
Notions of pseudorandomness are usually ``symmetric'' or ``global'' in their definitions as in the definition in \cite{thom1} or in the definition of an $(n,d,\lambda)$ graph. This latter notion is at first glance somewhat asymmetric in the sense that the conditions imposed on the degrees and codegrees are only for the vertices of $X$.  However, it is a simple exercise (which we shall not get into here) to show that these conditions also imply certain restrictions on the degrees and codegrees of the vertices of $Y$ as a consequence of the following analogue of the expander-mixing lemma (restricted to our setup):
\begin{theorem}[Theorem 2 in \cite{thom2}]\label{thm:thom}
Let $G(X,Y)$ be a bipartite graph with $|X|= k\leq n = |Y|$, which is Thomason pseudorandom with parameters $(p,\varepsilon)$. Then for every subset $A \subseteq X$ of size at least $1/p$ and every subset $B \subseteq Y$, with $|A| = a$ and $|B| = b$,
\[
|e(A,B) - pab| \leq \sqrt{pnab(1+\varepsilon pa)}.
\]
%where $e(A,B)$ denotes the number of edges having one end point in $A$ and the other in $B$. 
\end{theorem}
Again, we remark that Thomason's theorem in \cite{thom2} is stated for pseudorandom bipartite graphs $G(X,Y)$ with $|X|=|Y|=n$ and parameters $(p,\mu)$. But a glance at the proof there immediately tells us that the same proof works in our general setup as well. The interesting point is that this asymmetric definition of pseudorandomness also yields the aforementioned  theorem. A heuristic and somewhat simplistic explanation  for this is that we are restricting ourselves to bipartite graphs, and it is precisely  due to   the bipartite structure of the graph that the arguments go through.

Another reason why we prefer to work with this notion of pseudorandomness is that it is combinatorial in its definition; it only considers the degrees of the vertices and codegrees of pairs of vertices of $X$, which is  computationally easy to verify. In addition, it is a reasonably \emph{robust} notion which also allows us to generate several non-trivial examples of Thomason pseudorandom graphs. While it is true that many notions of pseudorandomness do pass onto subgraphs, we did not find any concrete statement in the literature that established the same here for this notion. % of Thomason pseudorandomness. 
So we took it on ourselves to prove its robustness; see the lemma in the Appendix for a precise statement. 

Pseudorandom graphs enjoy several very interesting properties. It is not hard to show that $(n,d,\lambda)$ graphs with $d-\lambda\geq 2$ are $d$-edge connected and as a simple consequence, it follows that for even $n$, $(n,d,\lambda)$ graphs have a perfect matching \cite{KS}. In the more general context, it is conceivable that Thomason pseudorandom graphs admit ``almost-perfect'' matchings, i.e., admit a perfect matching on at least $(1-o(1))|V|$ vertices under not-too-restrictive conditions. The second result of our paper proves a more general version of this statement for NMP for Thomason pseudorandom graphs.

Before we formally state our result, we need the following definition.
\begin{definition}[NMP-Approximability]
Suppose $\ep>0$. For functions $f,g:\R^+\rightarrow\R^+$ such that $f(x),g(x)\to 0$ as $x\to 0$, a bipartite graph $G(X,Y)$  is said to be $(f,g,\varepsilon)$-NMP approximable if there are subsets $\X \subseteq X$ and $\Y\subseteq Y$ such that:\begin{itemize}
\item  $\frac{|\X|}{|X|} \leq f(\varepsilon)$, $\frac{|\Y|}{|Y|} \leq g(\varepsilon)$
\item  The bipartite subgraph induced on the sets $X\setminus \X$ and $Y\setminus \Y$ has NMP.\end{itemize}
\end{definition} 
We now state our second main result of the paper.
\begin{theorem}\label{thm:pseudo}
Suppose $0<\varepsilon <1$, and let $\omega:\N\to\R^+$ be a non-negative valued function that satisfies $\omega(k)\to\infty$ as $k\to\infty$.  There exists an integer $k_0 = k_0(\varepsilon,\omega)$ such that the following holds. Suppose $p\ge\frac{\omega(k)}{k}$, $|X|=k, |Y|=n$ with $k_0< k\leq n$, and suppose $G=G(X,Y)$ is a Thomason pseudorandom bipartite graph with parameters $(p,\varepsilon)$. Then $G$ is $(f,g,\varepsilon)$-NMP-approximable with %$\omega$ is any function such that $\omega(x) \longrightarrow \infty$ as $x\longrightarrow\infty$ and where
\begin{enumerate}%[label=(\alph*)]
\item[(a)]\label{part:one} $f(x) = O(x)$, $g(x) = O(\sqrt{x})$ if $n> \frac{k}{\sqrt{\varepsilon}}$ and %assuming $p_0\geq \frac{3}{\sqrt{k}}$ and $\varepsilon_0\leq \varepsilon$; 
\item[(b)]\label{part:two} $f(x) = g(x) = O(\sqrt[4]{x}\log\left( \frac{1}{x}\right))$ if $n\leq \frac{k}{\sqrt{\varepsilon}}$.% and assuming $p_0\geq \frac{\omega(k)}{\sqrt{k}}$ and $\varepsilon_0\leq \varepsilon^2$,  %$\omega(x) \longrightarrow \infty$ as $x\longrightarrow\infty$.
\end{enumerate}
\end{theorem}
Note that in the statement of Theorem  \ref{thm:pseudo} the bounds $f=g=O(x^{1/4}\log(1/x))$ work for all $(k,n)$. The first part of the theorem is a stronger conclusion when $n\gg k$. At the level of generality of the statement of Theorem \ref{thm:pseudo}, it may in fact be {\it necessary} to delete some vertices from the graph in order to achieve NMP. Indeed, the definition of a Thomason pseudorandom graph does not preclude the existence of isolated vertices; in fact, one could add a few isolated vertices to $Y$ to get another pseudorandom graph with only slightly worse parameters! Also, on a less frivolous note, suppose $n=O(k)$ and $\omega(k)\ll\log k$, and consider $\bG(k,n,p)$; a consequence of the proof of the second item of Theorem \ref{thresh_NMP} (which appears later in the paper as Lemma \ref{lem:ub}) shows that there are isolated vertices in $Y$ {\it whp}. Since $\bG(k,n,p)$ is also Thomason pseudorandom {\it whp} it follows that over the sparser regime for $p$ (where Theorem \ref{thm:pseudo} is applicable),  the deletion of some vertices is indeed necessary to arrive at the conclusion of Theorem \ref{thm:pseudo}.

Theorem \ref{thm:pseudo} essentially says that if we have a not-too-sparse pseudorandom bipartite graph, i.e., a Thomason pseudorandom graph with  $p$ not too small, then we can remove a small fraction of vertices from both parts such that the graph induced by the remaining vertices has the normalized matching property. The sense of how small these sets are is described using the notion of NMP-Approximability defined above. As we shall see, the proof actually establishes an ``approximate decomposition'' theorem: the vertex set of any Thomason pseudorandom bipartite graph almost admits a decomposition into copies of what we call a \emph{Euclidean Tree} - a small tree that arises canonically via the execution of the Euclidean algorithm. Furthermore, the entire process of obtaining $\X$ and $\Y$ is algorithmic (and efficient) in nature and we consider this to be a major feature of our argument. After the publishing of this article, we have learned that this notion of Euclidean Trees has been defined prior to our work in the context of graphic matroids\footnote{We thank Attila Sali for bringing this to our attention.} (see \cite{SS}). So we find it quite interesting to see it reappear in the context of a seemingly unrelated problem.
 
The rest of the paper is organised as follows. The next section gives some preliminaries and sets up terminology and tools that will be of use in the latter sections. In Section \ref{sec:threshold:NMP} we prove Theorem \ref{thresh_NMP}, and in Section \ref{pseudo:sec}, we prove Theorem \ref{thm:pseudo}. The paper concludes with some remarks and open questions in Section \ref{sec:conc}, and an Appendix. As mentioned earlier, %There were two reasons for including this. While it is true that many notions of pseudo-randomness do pass onto subgraphs, we did not nd any concrete statement anywhere that established the same here in this notion of Thomason pseudorandom. Secondly, 
the lemma in the Appendix can serve as a generator of several examples of Thomason-pseudorandom graphs for which
Theorem \ref{thm:pseudo} is applicable. The main reason for including the Lemma  is that most of
the standard and well-studied examples of pseudorandom graphs that arise from algebraic structures/posets tend to have $|X| = |Y|$, or even in the cases where $|X| \neq |Y|$, the corresponding
bipartite graphs are much sparser than the ones we need in our hypothesis. %But for the sake of convenience of the reader, we have moved its statement and proof to the Appendix.

\section{Preliminaries}\label{prelim}
%Before we start, we set up some notation.  
Suppose $G(X,Y, E)$ is a bipartite graph. For $U\subseteq X\cup Y$, set $U_X:=U\cap X$, $U_Y:=U\cap Y$. For sets $A\subseteq X, B\subseteq Y$, by $G(A,B)$ we shall mean the subgraph of $G$ induced by the vertex set $A\cup B$. For a vertex $x$, $d(x)$ shall denote its degree, and for sets $A\subseteq X, B\subseteq Y$, $e(A,B)$ shall denote the number of edges between $A$ and $B$.

We shall repeatedly make use of the Chernoff bound:
\begin{theorem}\label{thm:Chernoff} [Chernoff Bound] (As in \cite{rg}) Suppose $X\sim Bin(n,p)$ is a binomial random variable and $\lambda:=\bE(X)= np$. Then for $t>0$ 
\begin{eqnarray*} \bP(X\ge \bE(X) + t) &\le& \exp\left(-\frac{t^2}{2(\lambda+t/3)}\right)\\
  \bP(X\le \bE(X) - t) &\le& \exp\left(-\frac{t^2}{2\lambda}\right).\end{eqnarray*}\end{theorem}

A natural question that arises in the context  of NMP is: If $G(X,Y)$ has NMP, then does $G(Y,X)$ also have NMP, i.e., is it true that for all $T\subseteq Y,\  \frac{|N(T)|}{|X|}\ge\frac{|T|}{|Y|}$? This is not immediately obvious from the definition of NMP, but it is indeed the case, as can be immediately seen from the second characterization of Theorem \ref{thm:kleitman} which is symmetric in $X$ and $Y$.

We begin with a simple proposition that will be instrumental in our proof of Theorem \ref{thresh_NMP} in Section \ref{sec:threshold:NMP}. For a graph $G(X,Y)$ that does not have NMP we say that a set of vertices $S\subseteq X$ {\it witnesses the violation of NMP for $G(X,Y)$} if $\frac{|N(S)|}{|Y|}<\frac{|S|}{|X|}$.
\begin{lemma}\label{NMP:witness}
Suppose $G(X,Y)$ with $|X|=k$, $|Y|=n$ does not have NMP. Then, if $T\subset Y$ witnesses the violation of NMP for $G(Y,X)$, then $X\setminus N(T)\subset X$ witnesses the violation of NMP for $G(X,Y)$.
Moreover, either there exists $S\subset X$ that witnesses the violation of NMP for $G(X,Y)$ with $|S|\le \frac{k}{2}$, or there exists $T\subset Y$ that witnesses the violation of NMP for $G(Y,X)$ with $|T|<\frac{n}{2}+\frac{n}{k}$.
\end{lemma}
\begin{proof}
If $T\subset Y$ witnesses the violation of NMP for $G(Y,X)$, then \[\frac{|N(T)|}{|X|}<\frac{|T|}{|Y|}\Rightarrow \frac{|X\setminus N(T)|}{|X|}>\frac{|Y\setminus T|}{|Y|}\geq \frac{|N(X\setminus N(T))|}{|Y|},\] where we subtracted both sides from $1$ and used the simple fact that $N(X\setminus N(T))\subseteq Y\setminus T$ in the final inequality. Now, to see the ``moreover" part,
as $G$ does not have NMP,  first let $S$ be a minimal set that witnesses the violation of NMP for $G(X,Y)$.  By the minimality of $S$, we have $|N(S)|\ge\frac{n}{k}(|S|-1)$. If $|S|\le\frac{k}{2}$, then we are through, so suppose that $|S|>\frac{k}{2}$. Let $T=Y\setminus N(S)$. Then note that $|T|<\frac{n}{2}+\frac{n}{k}$. But then by the argument above (which is symmetric in $X$ and $Y$), $T$ witnesses the violation of NMP for $G(Y,X)$. 

\end{proof}

We also take note of a couple of facts from literature on random graphs that will be useful in the proof of Theorem \ref{thresh_NMP}. By $d(x)$ (respectively $d(y)$) we mean the degree of vertex $x$ into $Y$ (respectively the degree of vertex $y$ into $X$) in $G(X,Y) = \bG(k,n,p)$.
\begin{fact}\label{large_mindeg1} Let $p\geq \frac{(1+\varepsilon)\log n}{k}$. For any fixed $r\in\bN$, in $G(X,Y)$, $d(x)\ge r$ for all $x\in X$ and $d(y)\ge r$ for all $y\in Y$ {\it whp}. \end{fact}
This follows from the following well known result (see \cite{Bol} for instance, chapter 3) that in $\bG(n,n,p)$ if $p=\frac{\log n+(r-1)\log\log n+\omega(n)}{n}$ for any function $\omega(n)$ that goes to infinity with $n$, then {\it whp} $\bG(n,n,p)$ has minimum degree $r$ since the number of vertices of degree $r$ is approximately Poisson. The same argument extends to $\bG (k,n,p)$ as well. %The primary use of this fact (which we apply with $r = \frac{C}{\ep}$ for some constant $C$) in the proof of Theorem \ref{thresh_NMP} is to let us assume that {\it whp}, the minimal witnesses obtained from Lemma \ref{NMP:witness} are not too small i.e., have size at least $12/\ep$, say.  The following fact is used in conjunction with the above.

\begin{fact}\label{large_mindeg2} Let $p\geq \frac{(1+\varepsilon)\log n}{k}$ and suppose $n\ge 2k$. Then in $G(X,Y)$, {\it whp} every $x\in X$ has degree at least $\frac{\ep n\log n}{2k}$. \end{fact}
This is an easy consequence of the Chernoff bound (Theorem~\ref{thm:Chernoff}). Indeed, since $\bE[d(x)]=(1+\ep)\frac{n\log n}{k}$, it follows that 
$$\bP\left[d(x)<\frac{\ep n\log n}{2k}\textrm{\ for\ some\ }x\in X\right]\le k\exp\left(-\frac{(1+\ep/2)^2n\log n}{2(1+\ep)k}\right)\le n^{-\ep^2/8}.$$

We now introduce an important ingredient that is vital to the proof of Theorem \ref{thm:pseudo}. 
 Suppose $\ell, L$ are positive integers with $gcd(\ell, L)=1$.  A tree will be called a {\it left-right tree} if the two color classes of its vertex set are labelled as ``left'' and ``right'' respectively. Since a connected bipartite graph admits a unique $2$-coloring of its vertices, a left-right tree can be thought of a tree with a label on each vertex denoting its color class.

\begin{figure}[ht]
\centering
\includegraphics[height = 7cm]{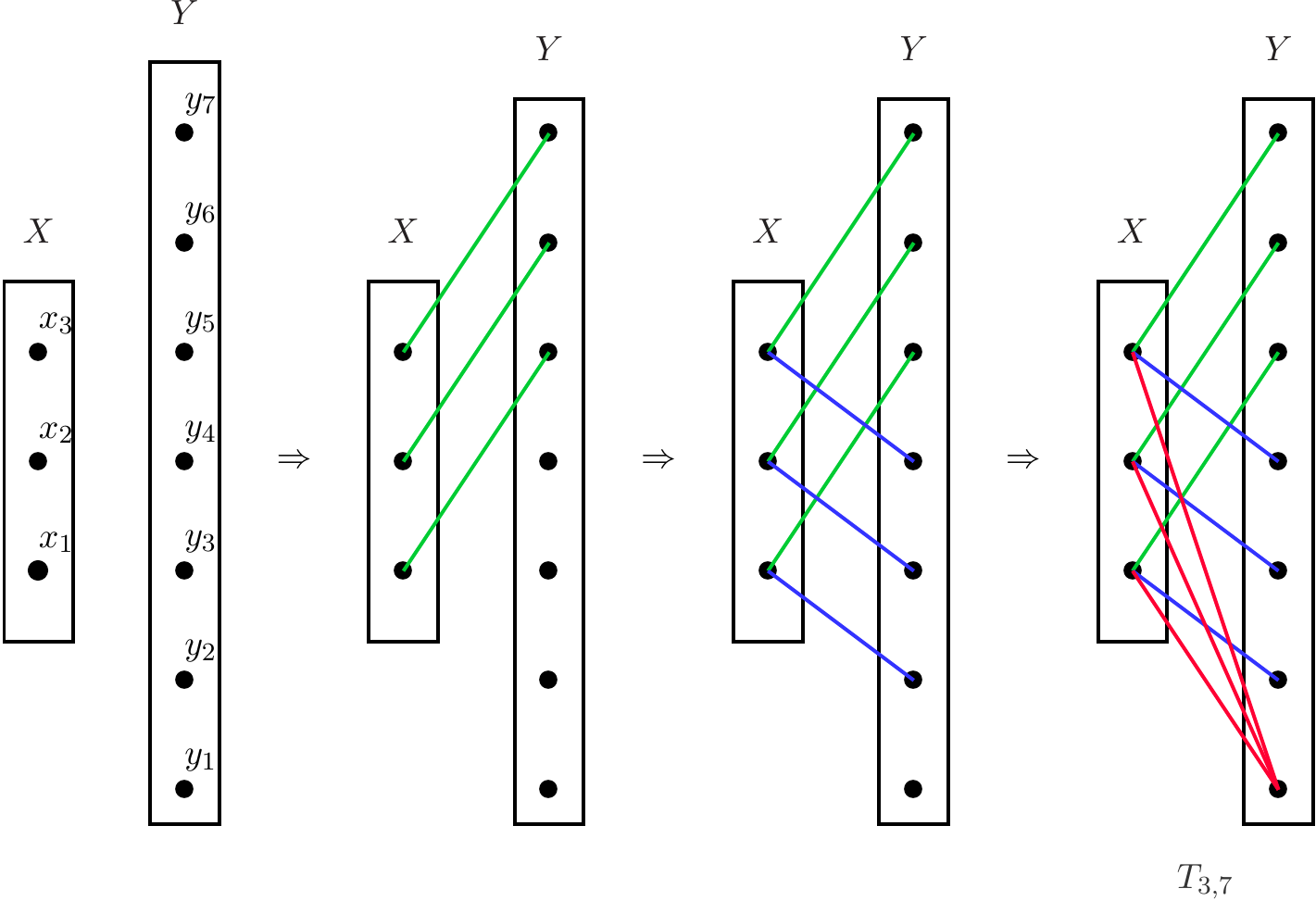}
\caption{Construction of the Euclidean $(3,7)$-tree. Each successive matching is shown in a different color.}\label{fig:37}
\end{figure}

\begin{figure}[ht]
\centering
\includegraphics[height = 7cm]{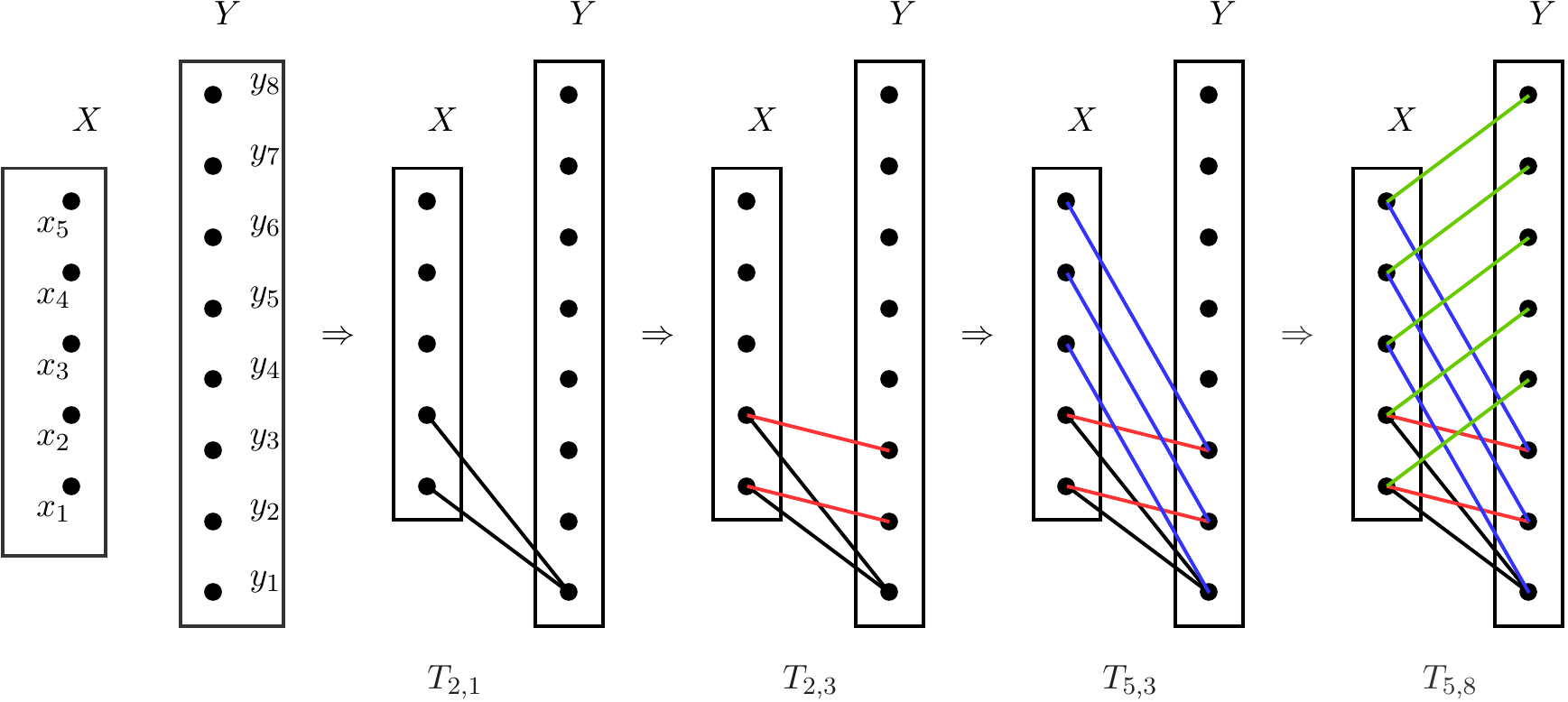}
\caption{The Euclidean $(5,8)$-tree process. In this case $m=4$, $(r_{2}, r_3, r_4,r_5)=(2,3,5,8)$, $(q_1,q_2,q_3,q_4)=(2,1,1,1)$.  $T_{5,8}$ evolves as $T_{2,1}\Rightarrow T_{2,3}\Rightarrow T_{5,3}\Rightarrow T_{5,8}$ in the process.}\label{fig:58}
\end{figure}

 The {\bf Euclidean $\mathbf{(\ell, L)}$-tree} which we shall denote by $T_{\ell, L}$, is a left-right tree on $\ell+L$ vertices with $\ell$ left vertices, and $L$ right vertices that is defined recursively as follows. If $\ell=1$, $T_{1,L}$ is simply a star on $L+1$ vertices with one left vertex and $L$ right vertices. If $L=1$, then $T_{\ell, 1}$ is the star on $\ell +1$ vertices with one right vertex, and $\ell$ left vertices. In general, suppose $X=\{x_1,\ldots,x_{\ell}\}$ and $Y=\{y_1,\ldots,y_L\}$ are the left and right vertex sets respectively, and suppose $\ell<L$.  Let $M_1$ denote the matching consisting of the edges $\{x_i, y_{i+L-\ell}\}$ for $1\le i\le \ell$. We define $T_{\ell, L}=M_1 \sqcup T_{\ell, L-\ell}$ where $\sqcup$ denotes an edge disjoint union, and $T_{\ell, L-\ell}$ is the corresponding Euclidean tree with left vertex set $X'=X$ and right vertex set $Y'=\{y_1,\ldots, y_{L-\ell}\}$. If $\ell >L$ then we define $M_1$ to be the matching $\{x_{i+\ell-L},y_i\}$ for all $1\le i\le L$ and define $T_{\ell,L}= M_1 \sqcup T_{\ell-L, L}$ where $T_{\ell-L,L}$ is the Euclidean tree with left vertex set $X'=\{x_1,\ldots,x_{\ell-L}\}$ and right vertex set $Y'=Y$.
A picture is worth a thousand words; see Figure \ref{fig:37} that illustrates the Euclidean tree $T_{3,7}$, and Figure \ref{fig:58} that illustrates $T_{5,8}$.

The following lemma conveys why Euclidean trees are relevant to us. 
\begin{lemma} \label{Euclidean:NMP} Suppose $T=T_{\ell, L}$ is a Euclidean tree. Then if $X, Y$ denote the sets of left and right vertices respectively, then $T$ as the bipartite graph $T(X,Y)$ has NMP. Moreover, so does the graph obtained by making several vertex-disjoint copies $T(X_i,Y_i)$ of $T$ i.e., the graph $\mathcal{T}(\mathcal{X},\mathcal{Y})$ where $\mathcal{X} = X_1\sqcup \cdots \sqcup X_r$, $\mathcal{Y} = Y_1\sqcup \cdots \sqcup Y_r$.
\end{lemma}
\begin{proof} First assume that $\ell<L$. If $\ell=1$, then $T$ is simply a star with $L$ leaves, and clearly, $T$ has NMP. Suppose by induction that Euclidean trees with fewer than $\ell+L$ vertices have NMP.  Let $S\subseteq X$. Then since  $T= M_1 \sqcup T_{\ell, L-\ell}$, it follows that $N(S)=\{y_{j+L-\ell}: x_j\in S\}\sqcup N'(S)$ where $N'(S)$ is the set of neighbors of $S$ among $\{y_1,\ldots,y_{\ell}\}$. But since $T_{\ell,L-\ell}$ has NMP, we have $|N'(S)|\ge \frac{L-\ell}{\ell}|S|$, so that $|N(S)|\ge |S|+\frac{L-\ell}{\ell} |S|=\frac{L}{\ell}|S|$ and that completes the proof. If $\ell > L$, then the above argument works with $\ell$ swapped with $L$ throughout and the fact that $T(X,Y)$ has NMP if and only if $T(Y,X)$ does. Finally, the observation that $\mathcal{T}(\mathcal{X},\mathcal{Y})$ has NMP follows immediately from the third (multiplicity function) characterization of NMP in Theorem \ref{thm:kleitman}.
\end{proof}

We now describe what we call the ``Euclidean ($\ell, L)$-tree process'' which details a realization of the graphs $T_{\ell, L}$ through a series of steps, which along with the corresponding terminology we build here will be relevant in Section \ref{pseudo:sec} in the proof of Theorem \ref{thm:pseudo}. This description also justifies why we call them Euclidean trees.
Suppose $\ell<L$. Consider the Euclidean algorithm on the pair $(\ell, L)$ as follows. 
\begin{eqnarray*} \label{Euclidean}
L&=&q_m\ell +r_{m-1},\hspace{2cm} 0<r_{m-1}<r_m=\ell,\\
\ell&=&q_{m-1}r_{m-1}+r_{m-2},\hspace{0.95cm} 0<r_{m-2}<r_{m-1},\\
\cdots&=&\cdots\\
r_{3}&=&q_{2}r_{2}+r_{1},\hspace{2.45cm} 0<r_{1}<r_2,\\
r_{2}&=&q_1r_{1},\hspace{3.27cm} r_{1}=1.\end{eqnarray*}
 If we set $r_{m+1}=L, r_m=\ell, r_0=0$, then we may write the equalities above as $r_{i+1}=q_{i}r_i+r_{i-1}$ for $1\le i\le m$. $m$ is referred to as the {\it complexity} of the Euclidean algorithm for the parameters $(\ell, L)$.
 The following fact is well-known (see for instance, \cite{knuth}, page 360).
 \begin{fact} \label{fact:euclid} The complexity of the Euclidean algorithm with input parameters $(\ell, L)$ is at most $2.078\log L+0.6723$.
 \end{fact}
 
We now describe $T_{\ell, L}$ as the evolution of an inductive sequence of trees through $m$ stages ($m$ as above), and in order to do that, we need some additional terminology. By an {\it $X$ $q$-fan}, we mean the tree $T_{1,q}$ and by a $Y$ $q$-fan, we mean $T_{q,1}$. By an $X$ {\it $q$-thrill\footnote{The collective noun for fans is a thrill, so the nomenclature seemed appropriate.}  of size $r$} we mean a union of $r$ vertex disjoint $X$ $q$-fans, and a $Y$ $q$-thrill is defined analogously. For a fixed graph $F$, an {\it $F$-factor} in a graph $G$ is a spanning subgraph of $G$ consisting of vertex disjoint copies of $F$.  As an example, an $X$ $q$-thrill admits a factoring by $X$ $q$-fans.

By  definition, $T_{\ell, L}$ is inductively obtained through a sequence of edge disjoint unions of matchings, until we finally terminate in a tree $T_{q,1}$ or $T_{1,q}$, % as the case may be, 
for some $q$. We now invert this process. 

Suppose $m$ as described above in the Euclidean algorithm is even (the odd case is analogous). Let $T_1:= T_{r_{\scriptscriptstyle 2},r_{\scriptscriptstyle 1}}=T_{r_2,1}$.  Having inductively defined $T_{i-1}$ with left set $X^{(i-1)}$, right set $Y^{(i-1)}$ and edge set $E_{i-1}$,  we define $T_{i}$ as follows. If $i$ is even, then the vertex set of $T_{i}$ has left set $X^{(i)}:=\{x_1,\ldots, x_{r_{i}}\}$, right set $Y^{(i)}=\{y_1,\ldots,y_{r_{i+1}}\}$, and the edges of $T_{i}$ consist of the edges of $T_{i-1}$ along with an additional $X$ $q_{i}$-thrill of size $r_{i}$ between the vertices of $X^{(i-1)}$ and the vertices of $Y^{(i)}\setminus Y^{(i-1)}$.  If $i$ is odd, then $T_{i}$ has  left vertex set $X^{(i)}:=\{x_1,\ldots,x_{r_{i+1}}\}$, right vertex set $Y^{(i)}:=\{y_1,\ldots, y_{r_{i}}\}$ and the edges of $T_{i}$ consist of the edges of $T_{i-1}$ along with an additional $Y$ $q_{i}$-thrill of size $r_{i}$ between the vertices of $X^{(i)}\setminus X^{(i-1)}$ and the vertices of $Y^{(i-1)}$.  In simpler terms, it is the same construction but with the roles of the left and right sets reversed as per the parity of $i$.  The main point is that the graphs $T_{i}$ are precisely the Euclidean trees $T_{r_{(i+1)}, r_{i}}$ (or $T_{r_{i},r_{(i+1)}}$  depending on the parity of $i$) along with isolated vertices. While the inductive definition of the Euclidean tree $T_{\ell, L}$ appends one additional matching at each step, the Euclidean tree process accelerates this by adding a $q$-thrill for an appropriate $q$. In particular, $T_m$ is precisely $T_{\ell, L}$ and as we shall see in Section \ref{pseudo:sec}, it is particularly handy to think of $T_{\ell, L}$ as the end result of this evolving process. Figure \ref{fig:58} gives an illustration of this evolution for the Euclidean tree $T_{5,8}$.

\section{Threshold for NMP for \texorpdfstring{$\mathbb{G}(k,n,p)$}{G(k,n,p)}}
\label{sec:threshold:NMP}
In this section we prove Theorem \ref{thresh_NMP}, restated below for convenience.  Throughout this section, we shall write $\bG$ to denote $\bG(k,n,p)$.  Unless stated otherwise, we shall assume $k\le n\le \exp(o(k))$.
 %{\it (Statement of Theorem \ref{thresh_NMP})}:\\ 
\paragraph*{Theorem $\mathbf{1.2}$.}
{\it
Suppose $k\leq n(k) \le \exp(o(k))$, and let $0<\ep,\delta <1$. There exists $k_0 = k_0 (\ep,\delta)$ such that for $k\ge k_0(\ep,\delta)$
\begin{enumerate}
\item If $p\geq \frac{(1+\ep)\log n}{k}$ then $\bP[\mathbb{G}(k,n,p)\textrm{\  has\  NMP}]\ge 1-\delta$. 
\item If $p\leq \frac{(1-\ep)\log n}{k}$ then $\bP[\mathbb{G}(k,n,p)\textrm{\  has\  NMP}]\le \delta$.
\end{enumerate}
}

We establish item 2 first i.e., that if $p$ is below the threshold then {\it whp}, $\bG$ does not have NMP. The proof is straightforward as it simply shows  the existence of an isolated vertex in $Y$ {\it whp}. 

\begin{lemma}\label{lem:ub}
Suppose $n = n(k)$ %grows to infinity with $k$ and 
be such that $k\leq n(k)$ for all $k\in \mathbb{N}$. Let $0< \ep < 1$. There exists $k_0 = k_0 (\ep)$ such that for $k\ge k_0$,
if $p\leq \frac{(1-\ep)\log n}{k}$ then $\mathbb{G}(k,n,p)$ does not have NMP {\it whp}.
\end{lemma}
%To complete the proof of Theorem \ref{thresh_NMP}, we shall show that if $p\leq \frac{(1-\ep)\log n}{k}$ then with high probability, $Y$ has an isolated vertex. 
\begin{proof}
Let $G(X,Y) = \bG$  and let $N$  denote the number of isolated vertices in $Y$. Then $\E[N] = n(1-p)^k$. 
\begin{claim}
Given $c>1$, there exists a unique $x_c\in (0,1)$ such that for all $x\in (0,x_c]$, $1-x\geq \exp(-cx)$ and equality holds only when $x = x_c$. Moreover, as $c\rightarrow 1^+$, $x_c\rightarrow 0^+$.
\end{claim}
The claim is a standard exercise in basic calculus, so we omit its proof. 
 
Fix $c$ such that $1<c<\frac{1}{1-\ep}$. Since $p < \frac{(1-\ep) \log n}{k} = o(1)$, by the above claim, there exists $k$ sufficiently large such that $1-p \geq \exp(-cp)$. Consequently, $\E[N] = n\cdot (1-p)^k \geq \exp(-cpk + \log n) = \exp(\alpha\log n) = n^\alpha$ which grows to infinity as $k$ does, where $\alpha = \alpha (\ep)$ is defined to be $1-c(1-\ep) > 0$. Now using the Chernoff bound (taking $t = \lambda = \E[N]$ in the second inequality in Theorem~\ref{thm:Chernoff}),  we have
%We shall use the following version of the Chernoff Bound (see \cite{AS15}):
%\begin{lemma}[Chernoff bound]\label{lem:chernoff}
%Suppose $T_1,\ldots, T_n$ are independent random variables taking values in $\{0, 1\}$. Let $T$ denote their sum and let $\mu = \E[T]$ denote its expected value. Then for any $0<\ep'<1$, \[\Pr[T\leq (1-\ep')\mu ] \leq \exp\left(\frac{-\ep'^2\mu}{2}\right).\]
%\end{lemma}

\[
\Pr[N= 0]  \leq \exp \left(-\frac{\E[N]}{2}\right) \le \exp\left(-\frac{n^\alpha}{2}\right) = \exp(-n^{\Omega_\ep(1)})= o(1)\]
for large $n$. This concludes the proof.
\end{proof}

Lemma \ref{lem:ub} establishes that the right threshold for having NMP in $\bG$ must be at least as large as $\frac{\log n}{k}$. The following is a heuristic argument that suggests that it is \emph{exactly} $\frac{\log n}{k}$. As mentioned in the Introduction,  a classical result of Erd\H{o}s-R\'enyi states that a sharp threshold for the existence of a perfect matching in a bipartite graph $\mathbb{G}(n,n)$ is $p=\frac{\log n}{n}$. In our present situation, suppose $k$ divides $n$. Replicate each vertex of $X$ by a factor of $n/k$ to obtain the set $X'$. Define the graph $G'(X',Y)$ as follows. If $x'\in X'$ arises from the replication of the vertex $x\in X$,  then $x'y\in E(G')$ if and only if $xy\in E(G)$. It is a straightforward exercise to see that the original graph $G(X,Y)$ has NMP if and only if $G'(X',Y)$ satisfies Halls' condition, or equivalently, $G$ has NMP if and only if $G'$ has a perfect matching. If this new bipartite graph behaves likes $\bG(n,n,p)$ (which it isn't), then we need $p\sim\frac{\log n}{n}$ for the existence of a perfect matching. But since each vertex of $X$ has been blown up to $n/k$ copies, it is intuitive to expect that each vertex of $G$ behaves like the union of all these $n/k$ vertices bundled together, which suggests a threshold of $\frac{n}{k}\cdot\frac{\log n}{n}=\frac{\log n}{k}$. While this argument is just a heuristic, it suggests what the correct threshold ought to be, as we next show is indeed the case by establishing the remaining (and main) item $1$ of Theorem \ref{thresh_NMP}.

Here is an overview of the proof. Lemma \ref{lem:ratio-large} proves the theorem when $n/k$ is large (i.e., grows to infinity with $k$), and this part of the proof only takes recourse to Theorem \ref{thm:kleitman}. The general case  however is a little more delicate. The basic idea in the general case of the proof considers estimating the probability that there is a {\it minimal} set $S$ that violates  the NMP condition. In that sense,  our strategy follows a line of argument \'a la Erd\H{o}s-R\'enyi  but we need  some additional ideas and more careful analysis to carry it through to fruition.

%In what follows we shall assume that $k,n$ are both sufficiently large, and $k< n < \exp(k)$.

\begin{lemma}\label{lem:ratio-large}
Suppose $n = k\omega(k)$ where the function $\omega(k)\geq 1$ for all $k\in \mathbb{N}$ and satisfies $\omega(k) \rightarrow \infty$ as $k\rightarrow\infty$. Let $0< \ep,\delta <1$. Then there exists $k_0 = k_0 (\ep,\delta)$ such that for $k\ge k_0(\ep,\delta)$, if $p\geq \frac{(1+\ep)\log n}{k}$, then $\bP[\mathbb{G}(k,n,p)\textrm{\  has\  NMP}]\ge 1-\delta$.
\end{lemma}

\begin{proof} Let $0 < \ep < 1/5$, and let $|X| = k\leq n  = |Y|$. Since NMP is a monotone property, it suffices to establish the lemma for $p=\frac{(1+\ep)\log n}{k}$.

%Let $\frac{n}{k}\ge\log n$. 
Suppose $\bG$ fails to have NMP. By Theorem \ref{thm:kleitman}, there exists an independent set $I=I_X\cup I_Y$ in $\bG$ such that
\( \frac{|I_X|}{k} + \frac{|I_Y|}{n} > 1. \)
Thus, from the union bound, the probability that $\bG$ does not have NMP is at most $\sum_{\ell = 1}^k P_\ell$ where for $1\leq \ell\le k$, where 
\begin{eqnarray}
P_\ell &=& \binom{k}{\ell}  \binom{ n}{ \left\lceil n\left(1 -\frac{\ell}{k}\right)\right\rceil}  (1-p)^{\ell \left\lceil n\left(1 -\frac{\ell}{k}\right)\right\rceil}\textrm{\ for\ }\ell<k\label{eq:P_ell}\\
P_k &=&n\cdot (1-p)^k \leq \exp(-(1+\ep)\log n + \log n) \leq \frac{1}{n^{\ep}}.\end{eqnarray}
Here, $P_\ell$ is an upper bound on the probability that there is a set $S\subseteq X$ of size $\ell$ and a set $T\subseteq Y$ of size $\left\lceil n\left(1 -\frac{\ell}{k}\right)\right\rceil$ such that $S\cup T$ is an independent set. $P_k$ is an upper bound on the probability that $Y$ contains an isolated vertex. 

We define $\ep':= \ep/2$ and split $\sum_{\ell} P_{\ell}$ into three cases according to whether $\ell$ is ``small'', ``intermediate'', or ``large'' and repeatedly make use of the well-known bounds $1+x\leq \exp (x)$ for all $x\in \mathbb{R}$ and the binomial coefficients $\binom{N}{K} \leq\left( \frac{eN}{K}\right)^K$  for all $K\le N$.
\paragraph{Small Case: $1\leq \ell \leq \ep' k$.} Here, using $\binom{ n}{ \left\lceil n\left(1 -\frac{\ell}{k}\right)\right\rceil} = \binom{ n}{ \left\lfloor \frac{n\ell}{k}\right\rfloor}$ followed by standard binomial coefficient bounds, (\ref{eq:P_ell}) yields
% Using the bounds $\binom{k}{\ell} \leq (\frac{ke}{\ell})^\ell$, $\binom{ n}{ \lceil n\left(1 -\frac{\ell}{k}\right)\rceil} = \binom{ n}{ \lfloor \frac{n\ell}{k}\rfloor} \leq (\frac{2ke}{\ell})^{\frac{n\ell}{k}}$, and  $1-p \leq \exp(-p)$ to get
\begin{align}
P_\ell &\leq \exp\left(\frac{-(1+\ep)\ell\log n}{k}\cdot \left\lceil n\left(1 -\frac{\ell}{k}\right)\right\rceil + \ell\cdot\left(1+\frac{n}{k}\right)\cdot\left(1 +\log \frac{k}{\ell}\right) \right)\\
&\leq \exp\left(-(1+\ep)\cdot n\log n \cdot \frac{\ell}{k}\left(1 -\frac{\ell}{k}\right) + \left(1+\frac{\ep}{8}\right)^2\cdot n\cdot \frac{\ell}{k}\cdot\log k \right)\label{eq:small1}\\
&\le \exp\left( \frac{n\ell}{k}\cdot \log n \left[-(1+\ep) (1-\ep') + \left(1+\frac{\ep}{8}\right)^2\right]\right)\label{eq:small2}\\
&< \exp\left(-\frac{\ep}{8}\cdot \frac{n}{k}\cdot\log n\right)\label{eq:small3}
\end{align}
where to derive (\ref{eq:small1}), we use the bounds $\left\lceil n\left(1 -\frac{\ell}{k}\right)\right\rceil \geq n\left(1 -\frac{\ell}{k}\right)$, $1+\log \frac{k}{\ell} \leq 1 + \log k\leq (1+\frac{\ep}{8})\log k$ and $1+\frac{n}{k} \leq (1+\frac{\ep}{8}) \frac{n}{k}$ for large enough $k$. This is where we crucially use our assumption that $n/k \rightarrow \infty$ as $k\rightarrow\infty$.  (\ref{eq:small2}) follows by using the trivial fact that  $\log k \leq \log n$ and taking out the common factor $\frac{n\ell}{k}\cdot \log n$. 
%\[
%P_\ell \leq .
%\]
% Upon plugging in $\ep' = \ep/2$ and $\ell \geq 1$, we see that 
%\[P_\ell \leq \exp\left( \left[\frac{33\ep^2}{64}- \frac{\ep}{4}\right]\cdot\frac{n}{k}\cdot\log n\right)< \exp\left(-\frac{47\ep}{320}\cdot \frac{n}{k}\cdot\log n\right)<\] where we have used the assumption that $\ep<1/5$ in the middle inequality.
(\ref{eq:small3}) is obtained by using $\ell \geq 1$, plugging in $\ep' = \ep/2$ and working out that the expression in the square brackets in (\ref{eq:small2}) is at most $-\ep/8$ for small $\ep$. Finally, since $\frac{n}{k} > \frac{16}{\ep}$  for large enough $k$, it follows that $P_\ell < 1/n^2$ in this case. 
%Next, we also show that for this setting of $\ep'$, the combined error of the remaining cases on $\ell$ is also at most inverse polynomial in $n$, thereby finishing the proof.
\paragraph{Intermediate Case: $\ep' k \leq \ell \leq (1-\ep') k$.} Using the same expression for the upper bound on $P_\ell$ as in the previous case, we have
\[P_\ell \leq \exp\left(-(1+\ep)\cdot n\log n \cdot \frac{\ell}{k}\left(1 -\frac{\ell}{k}\right) + \ell\cdot\left(1+\frac{n}{k}\right)\cdot\left(1 +\log \frac{k}{\ell}\right) \right)
\]
Using the observation that in this case, $\frac{\ell}{k}(1-\frac{\ell}{k}) \geq \ep'(1-\ep')$ and the trivial bound $1+\frac{n}{k}\leq \frac{2n}{k}$, we obtain
\[ P_\ell \leq \exp\left(-(1+\ep)\cdot n\log n \cdot \ep'(1-\ep') + 2 n\cdot \frac{\ell}{k}\cdot\log \frac{k}{\ell} + \frac{2n\ell}{k} \right)< \exp\left(-\frac{\ep n\log n }{2}  + 3 n \right)\]
where the last inequality follows - setting $x = \ell/k$ - from the fact that $x\log \frac{1}{x} < 0.5$ for all $0<x<1$. Hence, $P_\ell < \frac{1}{n^{\ep n/3}}.$% implies that $P_\ell \leq \frac{1}{n^{\ep n}}$ for $n$ sufficiently large. 

\paragraph{Large Case: $(1-\ep') k \leq \ell < k$.} This case is completely analogous to the small case.  %we start with equation~\ref{eqn:bd} and bound each of the terms on the right hand side. 
First, observe
\[n\left(1 -\frac{\ell}{k}\right)\leq\left\lceil n\left(1 -\frac{\ell}{k}\right)\right\rceil\leq \left(1+\frac{\ep}{8}\right)n\left(1 -\frac{\ell}{k}\right)\] for large enough $k$ (again using $n/k \rightarrow \infty$ as $k\rightarrow\infty$) and we have that $P_\ell$ is at most
%where the last inequality is true for large enough $n$ and $k$ as $n = \omega(k)$. 
%We write down the following bound on the binomial coefficient 
%$$\binom{ n}{ \lceil n\left(1 -\frac{\ell}{k}\right)\rceil}  \leq (\frac{ke}{k -\ell})^{\lceil n\left(1 -\frac{\ell}{k}\right)\rceil}\leq (\frac{ke}{k -\ell})^{\left(1+\frac{\ep}{8}\right)n\left(1 -\frac{\ell}{k}\right)}$$
%$\binom{ n}{ \lceil n\left(1 -\frac{\ell}{k}\right)\rceil}  \leq (\frac{ke}{k -\ell})^{ n\left(1 -\frac{\ell}{k}\right)}$
% and that $\binom{k}{\ell} = \binom{k}{k-\ell} \leq (\frac{ke}{k - \ell})^{k-\ell}$. Plugging these bounds in equation~\ref{eqn:bd}, we obtain 
\begin{align*}
&\exp\left(-(1+\ep)\cdot \log n \cdot \frac{\ell}{k}\cdot \left\lceil n\left(1 -\frac{\ell}{k}\right)\right\rceil +(k-\ell)\left(1+\left(1+\frac{\ep}{8}\right)\frac{n}{k}\right)\left(1+\log \frac{k}{k-\ell}\right) \right)\\
&\leq \exp\left(-(1+\ep)\cdot n\log n \cdot \frac{\ell}{k}\left(1 -\frac{\ell}{k}\right) + \left(1+\frac{\ep}{8}\right)^3n\cdot\left(1- \frac{\ell}{k}\right)\cdot\log k \right)
\end{align*}
where in the last step we use the bound $1+\log \frac{k}{k - \ell} \leq 1+ \log k\leq \left(1+\frac{\ep}{8}\right)\log k$ for large enough $k$. Consequently, 
%we use $\log k \leq \log n$ and pull out the common factor $n \log n\cdot \left(1 -\frac{\ell}{k}\right)$ to obtain
\begin{eqnarray*}
P_\ell &\leq& \exp\left( n \log n\cdot \left(1 -\frac{\ell}{k}\right) \left[-(1+\ep) (1-\ep') + \left(1+\frac{\ep}{8}\right)^3\right]\right)\\
&\le& \exp\left( \frac{n \log n}{k}\cdot \left[-(1+\ep) (1-\ep') + \left(1+\frac{\ep}{8}\right)^3\right]\right) \le \frac{1}{n^2}.
\end{eqnarray*}
To explain the last step, the expression within the square brackets evaluates to $\frac{\ep}{512}(\ep^2+280\ep-64)$ which is at most $\frac{-199\ep}{12800}< \frac{-\ep}{128}$ when $0<\ep<1/5$. But $\frac{n}{k}>256/\ep$ for  sufficiently large $k$ and $n$ since $n/k\rightarrow\infty$. Thus, we have $\sum_{\ell} P_{\ell}=o(1)$ and that completes the proof of the lemma.
\end{proof}

Note that the argument in the intermediate case does not require $k = o(n)$ and in fact shows the following (in light of Theorem \ref{thm:kleitman}, switching from the independent set viewpoint to the violation of NMP viewpoint):
\begin{corollary}\label{remark:caseII} Given $\ep > 0$, for any $k\le n$ large enough, and vertex sets $X$ and $Y$ of sizes $k$ and $n$ respectively, the probability that there exists $S\subset X$ with $\ep' k\le |S|\le (1-\ep')k$ for $\ep' = \ep/2$ such that $S$ witnesses a violation of NMP for $G(X,Y) = \bG (k,n,p)$ is at most $n^{-\Omega_\ep (n)}$.  \end{corollary}

Interestingly, the proof of Lemma \ref{lem:ratio-large} actually works out \emph{for all} $n\geq k$ if one assumes $p\geq \frac{10\log n}{k}$ in the hypothesis instead of the sharper assumption on $p$. This, combined with Lemma \ref{lem:ub}, already establishes that $\frac{\log n}{k}$ is a \emph{threshold} for NMP. The additional ideas employed in the remainder of this section are essentially only required to show that $\frac{\log n}{k}$ is a \emph{sharp} threshold.

\begin{proof}[Proof of Theorem \ref{thresh_NMP}.]
In light of Lemma \ref{lem:ratio-large}, it suffices to prove the theorem assuming $\frac{n}{k}\le \log n$. $\log n$ here may be replaced by any slow-growing (but unbounded) function of $k$ or $n$ without much change to the rest of the argument, but we stick to $\log n$ for convenience.  %We start with a couple of facts. 
 
By Lemma \ref{NMP:witness} either there exists $S\subset X$ with $|S|\le k/2$ that witnesses a violation of NMP for $G(X,Y)$, or there exists $T\subset Y$ with $|T|<\frac{n}{2}+\frac{n}{k}$ that witnesses the violation of NMP for $G(Y,X)$ (of course, these cases need not be mutually exclusive; we merely use that combined, they exhaust the event that NMP is violated). The proof naturally splits into cases (labelled $X$ and $Y$ respectively)  according to whether the set winessing the violation is a subset of $X$ or $Y$. We shall show that either case occurs with low probability by exploiting certain properties of the minimal witness.  

\paragraph{Case $X$:} Define $\ell_{\min}$ to be the constant $\frac{18}{\varepsilon}$ if $1\leq \frac{n}{k} < 2$ and $\frac{\varepsilon \log n}{2}$ if $2\leq \frac{n}{k} \leq \log n$. In light of Facts~\ref{large_mindeg1} (for $r=\frac{36}{\ep}$ if $1\leq\frac{n}{k}<2$) and~\ref{large_mindeg2}, it follows that any minimal $S\subset X$ that witnesses the violation of NMP for $G(X,Y)$ must have size at least $|S| \geq \frac{k\delta(G)}{n}\geq \ell_{\min}$ \emph{whp} where $\delta(G)$ denotes the minimum degree of the vertices in $X$. The choice of the peculiar constant $r = \frac{36}{\varepsilon}$ will become clear later.

 Suppose $S\subset X$ such that $\ell_{\min}\leq |S|=\ell \leq \ep' k$ where $\ep' = \frac{\ep}{2}$. We first claim that every $U\subset N(S)$ of size $\left\lceil\frac{n}{k}\right\rceil$ witnesses at least $2$ neighbors (as a set) in $S$.  Indeed, suppose there is a subset $U$ of $\lceil\frac{n}{k}\rceil$ vertices in $N(S)$ which are the neighbors of only one vertex $x$ in $S$. Then by the minimality of $S$, it follows that the set $S'=S\setminus\{x\}$ satisfies $\frac{n}{k}|S|-\lceil\frac{n}{k}\rceil>|N(S')|\ge\frac{n}{k}(|S|-1)$ which is a contradiction, and that proves the claim. 

We divide case $X$ further into two subcases. First, we bound the probability that there exists $S\subset X$ of size $\ell$ for which $\frac{4\ell n\log n}{k^2}<1$ (notice that this clearly implies $\ell \leq \ep' k$) which witnesses a violation of NMP for $G(X,Y)$. So fix a choice for $S\subset X$ of size $\ell$, and $T\subset Y$ (which will represent $N(S)$) of size equal to some integer in the interval $[\frac{n\ell}{k}-\frac{n}{k},\frac{n\ell}{k})$. Fix a partition of $T$ into sets of size $\left\lceil\frac{n}{k}\right\rceil$. By size considerations, there are at least $t=\left\lfloor\frac{n(\ell-1)}{k\lceil n/k\rceil}\right\rfloor\ge\left\lfloor\frac{\ell-1}{1+(k/n)}\right\rfloor\ge\left\lfloor\frac{\ell-1}{2}\right\rfloor$ such parts, and by the observation above, each such part admits at least two neighbors in $S$. We conclude that the probability that there exists $S\subset X$ with $|S|\le\frac{k^2}{4n\log n}$ which witnesses a violation of NMP for $G(X,Y)$  is at most 

\begin{eqnarray}\label{probab_est}\Sigma_1 := \frac{n}{k}\sum_{\ell\geq \ell_{\min}} \binom{k}{\ell}\binom{n}{\lfloor\frac{n\ell}{k}\rfloor}(1-p)^{\ell \left\lceil n\left(1 -\frac{\ell}{k}\right)\right\rceil}\left(\binom{\ell}{2}\left(\left\lceil\frac{n}{k}\right\rceil p\right)^2\right)^t.\end{eqnarray} %where $S_1$ denotes the sum over $\ell$ for which $\frac{4\ell n\log n}{k^2}<1$, and $S_2$ ranges over the remaining $\ell$. 

To see why, observe that there are $\binom{k}{\ell}$ choices for $S$, at most $n/k$ values for $|N(S)|$ (since $S$ minimally witnesses a violation of NMP), each of which is at most $\lfloor\frac{n\ell}{k}\rfloor$. The probability that $e(S,Y\setminus N(S))=0$ is at most $(1-p)^{\ell \left\lceil n\left(1 -\frac{\ell}{k}\right)\right\rceil}$, and finally, the last expression is a bound on the probability that each of the $t$ blocks of vertices has at least $2$ neighbors in $S$. The condition on $\ell$ that we have imposed in this subcase simply translates to the observation that the quantity in the right-most parenthesis that is raised to $t$ is less than $1$. So, we have
%Let $r=\frac{18}{\ep}$.  %Then, by observation \ref{large_minddeg} $\bG$ has minimum degree at least $r$ with 
\begin{eqnarray*} 
\Sigma_1&\le&\frac{n}{k} \sum_{\ell\geq \ell_{\min}} \frac{(\frac{ek}{\ell})^{\ell}(\frac{ek}{\ell})^{(n\ell /k)}\left(\ell\lceil n/k\rceil p\right)^{2t}}{n^{(1+\ep)(n\ell/k)(1-\frac{\ell}{k})}}\\
\left(\text{Using }2t \geq \ell -3\text{ and }p\leq \frac{2\log n}{k}\right)      &\le& \frac{n}{k}\left(\frac{k^2}{4n\log n}\right)^3\sum_{\ell} \left(\frac{(\frac{ek}{\ell})^{(n/k)}\cdot\frac{ek}{\ell}\cdot\frac{4\ell n\log n}{k^2}}{n^{(1+\frac{\ep}{3})(n/k)} }\right)^{\ell}\\
%+(\frac{k^2}{4n\log n})^2\sum_{\ell\ge\frac{\ep\log n}{2}} \left(\frac{(\frac{ek}{\ell})^{n/k}\cdot\frac{ek}{\ell}\cdot(\frac{4\ell n\log n}{k^2})^2}{n^{(1+\frac{\ep}{3})n/k} }\right)^{\ell}
\left(\text{Using }\frac{n}{k}\leq \log n\right)     &\le & \frac{k^3}{64\log^3 n}\sum_{\ell\geq \ell_{\min}}  \left[\left(\frac{ek}{\ell n}\right)^{(n/k)}\cdot\left(\frac{4e\log^2 n}{n^{\ep/3}}\right)\right]^{\ell}\\
     &\le& \frac{k^3}{32\log^3 n}\left(\frac{4e\log^2 n}{n^{\ep/3}}\right)^{\ell_{\min}}
\end{eqnarray*}
for $n, k$ sufficiently large and where in the final step, we used the fact that an infinite geometric series is at most twice the first term, when the common ratio is small enough.
This expression is clearly $o(1)$ when $\frac{n}{k}\geq 2$ (and so $\ell_{\min} = \frac{\varepsilon\log n}{2}$). Further, it is at most $\frac{k^3}{32\log^3 n}\left(\frac{4e\log^2 n}{n^{\ep/6}}\right)^{18/\ep}=O(\frac{1}{\log^3 n})=o(1)$ when $1\leq\frac{n}{k}< 2$.
%If $\frac{n}{k}\geq 2$ then by Fact \ref{large_mindeg2} any minimal set witnessing a violation of NMP must have size at least $\ell\ge\frac{\ep\log n}{2}$ {\it whp}. Hence the last expression is clearly $o(1)$ for $n, k$  sufficiently large. If $1\le\frac{n}{k}\le 2$, set $r=\frac{36}{\ep}$. Then by Fact \ref{large_mindeg1} {\it whp}  $\bG$ has minimum degree at least $r$, so in particular, a minimal set $S$ witnessing the violating on NMP has size at least $\frac{rk}{n}\ge \frac{18}{\ep}$. 

For the subcase $\frac{k^2}{4n\log n}\leq \ell \leq \ep' k$, %the bound  is similar to that in equation~\ref{probab_est}, except that we forgo the last parenthesis. In other words, 
we simply bound (which we shall call $\Sigma_2$) the probability of a minimal $S$ whose size is in this range by the probability that $S\cup \overline{N(S)}$ is independent and sum over the entire range of $\ell$ again. First, observe that in this subcase, 
\[
\frac{ek}{\ell} \leq \frac{4en\log n}{k}\leq 4e\log^2 n
\] and thus,
\begin{eqnarray*}
\Sigma_2 &\le& \sum_{\ell} \frac{(\frac{ek}{\ell})^{\ell}(\frac{ek}{\ell})^{(n\ell /k)}}{n^{(1+\ep)(n\ell /k)(1-\frac{\ell}{k})}}\\
&\le& \sum_{\ell} \left[\frac{(\frac{ek}{\ell})^{1+(n/k)}}{n^{(1+\ep/3)(n/k)}}\right]^\ell\\
&\le& \sum_{\ell\geq \frac{k^2}{4n\log n}} \left[\left(\frac{4e\log^2 n}{n^{1+\ep/6}}\right)^{n/k}\cdot\left(\frac{4e\log^2n}{n^{\ep/6}}\right)\right]^{\ell} = o(1)\\
%&\le&  \frac{k^2}{16\log^2 n}\left(\frac{16e\log^4 n}{kn^{\ep/6}}\right)^{(\ep\log n)/2}= o(1)
\end{eqnarray*}
as before and we are through. 

Finally, observe that the case $\ep' k\le |S|\le k/2$ follows immediately from Corollary \ref{remark:caseII}.

\paragraph{Case $Y$:} There is a minimal witness $T\subset Y$ with $|T|=s\le \frac{n}{2}+\frac{n}{k}$ that witnesses the violation of NMP for $G(Y,X)$. This time though, since $k\leq n$ it follows that $|N(T)|\leq\lfloor\frac{ks}{n}\rfloor$, and that for every $x\in N(T)$ there are at least $2$ neighbors in $T$. Now, define $s_{\min} \coloneqq \frac{12}{\ep}$. As earlier, by Fact \ref{large_mindeg1}, the minimal $T\subset Y$ that witnesses the violation of NMP for $G(Y,X)$ must have size at least $s_{\min}$ \emph{whp}. Again, we split this into two subcases: $s_{\min} \leq s\le\ep' n$ and $s\ge\ep' n$ where again $\ep' =\ep/2$.

Suppose $s_{\min} \leq s\le\ep' n$. Analogous to how we divided Case $X$ into two subcases, let us first assume that $s\le \frac{k}{2\log n}$ which in particular, lets us assume that $sp < 1$. Then the probability that such a witness exists of size in this range is at most 
\begin{eqnarray}
M_1 &=& \sum_{s_{\min}\leq s\leq \frac{k}{2\log n}} \binom{n}{s}\binom{k}{\lfloor\frac{ks}{n}\rfloor}(1-p)^{s(k-\left\lfloor\frac{ks}{n}\right\rfloor)}\left(\binom{s}{2}p^2\right)^{\left\lfloor\frac{ks}{n}\right\rfloor}\\
%\left(\text{Using }\left\lfloor\frac{ks}{n}\right\rfloor \geq \frac{ks}{n} - 1 \text{ in the exponent and } \left\lfloor\frac{ks}{n}\right\rfloor \geq \frac{ks}{2n}\text{ elsewhere}\right) 
&\leq & \frac{1}{s_{\min}^2p^2}\sum_{s}\left(\frac{en}{s}\right)^s \left(\frac{2en}{s}\right)^{(ks/n)}(1-p)^{sk(1-\frac{s}{n})}(s^2p^2)^{\frac{ks}{n}}\label{eq:first}\\
%\left(\text{Using }\left\lfloor\frac{ks}{n}\right\rfloor \geq \frac{ks}{n} - 1\right) 
&\le& \frac{144 k^2}{\ep^2\log^2 n}\sum_{s}\left[\frac{\left(\frac{2en}{s}\right)^{1+(k/n)}\left(\frac{2s^2\log^2 n}{k^2}\right)^{(k/n)}}{\exp\left(pk\left(1-\frac{s}{n}\right)\right)}  \right]^{s}\\
\left(\text{Using $s\leq \frac{k}{2\log n}$}\right)&\leq&  \frac{144k^2}{\ep^2\log^2 n}\sum_{s}\left[\frac{\left(\frac{2en}{s}\right)\left(\frac{4esn\log^2 n}{k^2}\right)^{(k/n)}}{\exp\left((1+\ep)\log n - \frac{(1+\ep)k}{2n}\right)}  \right]^{s}\\
\left(\text{Using $s\leq \frac{k}{2\log n}$ again}\right)
 &\leq&  \frac{144k^2}{\ep^2\log^2 n}\sum_{s}\left[\frac{2en\cdot \left(\frac{2en\log n}{k}\right)^{(k/n)}}{n^{1+(\ep/2)}}  \right]^{s}\\
\left(\text{As $\frac{n}{k}\leq \log n$ and $\frac{k}{n}\leq 1$}\right) &\leq&  \frac{144k^2}{\ep^2\log^2 n}\sum_{s\geq s_{\min}}\left[\frac{4e^2\log^2 n}{n^{\ep/2}}\right]^{s}\\
\text{(Geometric series bound)} &<&  \frac{144k^2}{\ep^2\log^2 n}\left(\frac{8e^2\log^2 n}{n^{\ep/2}}\right)^{12/\ep} = o(1)
%&\leq&  \frac{k^{2-\frac{\varepsilon}{2}}}{\log^2 n}\left(\frac{4e^2\log^3 n}{n^{\ep/3}}\right)^{(\ep\log n)/2} = o(1)
\end{eqnarray}
where to derive (\ref{eq:first}), we use $\lfloor\frac{ks}{n}\rfloor\geq \frac{ks}{n} - 1$ in the exponent and the more crude bound $\left\lfloor\frac{ks}{n}\right\rfloor \geq \frac{ks}{2n}$ elsewhere, which is applicable since by assumption, $\left\lfloor\frac{ks}{n}\right\rfloor \geq |N(T)| \geq s_{\min} > 1$. We also subsequently drop the range $\frac{k}{2\log n}\geq s\geq s_{\min}$ in the sum for convenience. Next, if $ \frac{k}{2\log n}\leq s\leq \ep' n$, then we simply bound the probability of there being a witness of size in this range by the probability that $T\cup \overline{N(T)}$ is an independent set (i.e. the final parenthesis in the expression for $M_1$ above is dropped) and sum over this range of $s$ again. The calculations (for the accordingly defined expression $M_2$) are very similar to that of $\Sigma_2$ in case $X$ and are omitted here.

Finally, if $|T|>\ep' n$, then note that $S=X\setminus N(T)$ has size $(1-\ep')k\ge |S|\ge\ep' k$, and by Lemma \ref{NMP:witness},  $S$ witnesses the violation of NMP for $G(X,Y)$ and is covered by Corollary \ref{remark:caseII}.
\end{proof}

\section{Normalized Matching Property in Pseudorandom Graphs}\label{pseudo:sec}
%In this section, we discuss the normalized matching property in pseudorandom graphs. 
In this section, we prove Theorem \ref{thm:pseudo} which is restated below for convenience. Suppose $0<p<1$ and $0< \ep < 1$.  Recall that a bipartite graph $G(X,Y)$ with $|X|= k\leq n = |Y|$ is called Thomason pseudorandom with parameters $(p,\varepsilon)$ if every vertex in $X$ has degree at least $pn$, and if every pair of vertices in $X$ have at most $p^2n(1 +\varepsilon)$ neighbors in common.

\paragraph*{Theorem \ref{thm:pseudo}.}
{\it
Suppose $0<\varepsilon <1$, and let $\omega:\N\to\R^+$ be a non-negative valued function that satisfies $\omega(k)\to\infty$ as $k\to\infty$.  There exists an integer $k_0 = k_0(\varepsilon,\omega)$ such that the following holds. Suppose $p\ge\frac{\omega(k)}{k}$, $|X|=k, |Y|=n$ with $k_0< k\leq n$, and suppose $G=G(X,Y)$ is a Thomason pseudorandom bipartite graph with parameters $(p,\varepsilon)$. Then $G$ is $(f,g,\varepsilon)$-NMP-approximable with %$\omega$ is any function such that $\omega(x) \longrightarrow \infty$ as $x\longrightarrow\infty$ and where
\begin{enumerate}%[label=(\alph*)]
\item[(a)] $f(x) = O(x)$, $g(x) = O(\sqrt{x})$ if $n> \frac{k}{\sqrt{\varepsilon}}$ and %assuming $p_0\geq \frac{3}{\sqrt{k}}$ and $\varepsilon_0\leq \varepsilon$; 
\item[(b)] $f(x) = g(x) = O(\sqrt[4]{x}\log\left( \frac{1}{x}\right))$ if $n\leq \frac{k}{\sqrt{\varepsilon}}$.% and assuming $p_0\geq \frac{\omega(k)}{\sqrt{k}}$ and $\varepsilon_0\leq \varepsilon^2$,  %$\omega(x) \longrightarrow \infty$ as $x\longrightarrow\infty$.
\end{enumerate}
}
 In what follows, $G=G(X,Y)$ is a Thomason pseudorandom graph with parameters $(p,\ep)$ where $\ep > 0$ and $p\geq \frac{\omega(k)}{k}$ where $\omega(k)$ denotes a function that satisfies $\omega(k)\to\infty$ as $k\to\infty$. As always, $|X|=k\le n=|Y|$, and $n, k$ are sufficiently large (depending on the choice of $\ep$ and $\omega$). As in the proof of Theorem \ref{thresh_NMP}, we split the task of proving NMP-approximability into two cases: the first, in which $n$ is significantly larger than $k$ and the second, in which the two are comparable. 
 
 Here is a brief overview of the proof. Suppose that 
\[
\frac{n}{k} = \frac{L}{\ell},
\]
where the latter is the \emph{representation in reduced form} i.e., $gcd(\ell,L) = 1$ and $\ell, L \in \N$. Our strategy of proof is to show that we can find small sets $D_X\subset X, D_Y\subset Y$ such that $G(X\setminus D_X, Y\setminus D_Y)$ admits a vertex decomposition into copies of the Euclidean tree $T_{\ell, L}$. Since $T_{\ell, L}$ has NMP by Lemma \ref {Euclidean:NMP}, this establishes the NMP-approximability of $G$. An essential ingredient in the proof of both cases is Lemma~\ref{lem:main} (which appears below) which basically states: If $G(X,Y,E)$ satisfies that for every subset $A \subseteq X$ of size at least $1/p$ and every subset $B \subseteq Y$,
we have $|e(A,B) - p|A||B|| \leq \sqrt{pn|A||B|(1+\varepsilon p|A|)}$, then all large enough subsets of $X,Y$ admit an almost partition into $X$-thrills or $Y$-thrills (as the case may be).  %Since Thomason pseudorandom graphs satisfy this condition (Theorem \ref{thm:thom}), the proof is complete.

The application of this lemma in the first case ($n/k$ large) is straightforward, but in the second case, it does not apply directly. The principal issue in the second case emanates from the possibility that in the reduced form $\ell, L$ are still large; for instance if $n, k$ are coprime, then $(\ell, L)=(k,n)$ and Lemma \ref{lem:main} does not apply. To circumvent this difficulty, we pre-process the graph, by deleting a small portion from both $X, Y$ to get $X', Y'$ so that the reduced form $(\ell, L)$ for $(|X'|, |Y'|)$ satisfies $\ell, L=O_{\ep}(1)$.  Lemma \ref{lem:main} then applies in a multi-step process that we describe in Lemma \ref{lemma:smallcase}. 

 %(which will usually be clear from the context). 

\begin{lemma}\label{lem:main} %Suppose $q\in \N$ satisfies $pn - q > 0$ and let $\ep > 0$. %Let $G=G(X,Y)$ be Thomason pseudorandom with parameters $(p,\ep)$ where $p\ge\frac{\omega(k)}{k}$. 
Let $\ep > 0$ and $q\in\bN$ be such that $q = \left\lfloor \frac{n}{k}\right\rfloor$ or $q = O_\ep(1)$. %be a non-negative integer. 
Suppose $G(X,Y,E)$ satisfies the conclusion of Theorem \ref{thm:thom}. Let $U\subseteq X$ and $V\subseteq Y$ and define $d_0 = 2\ep n$. %$d_0 = 2\varepsilon n$. 
Then there exist subsets $A \subseteq U,B\subseteq V$ such that if $|U| = u, |V| = v, |A| = a,$ and $|B| = b$, then
\begin{itemize}
\item if $v = q u$, then $G(U\setminus A, V\setminus B)$ is spanned by an $X$ $q$-thrill %{\color{red} Not quite, correct this} 
where $a \leq d_0/q$ and $b \leq d_0$; 
\item if $u = q v$, then $G(U\setminus A, V\setminus B)$ is spanned by a $Y$ $q$-thrill where $a \leq q d_0$ and $b \leq  d_0$.
\end{itemize}
\end{lemma}
%{\color{red} In preliminaries, define G(A,B) as being the induced subgraph on A and B. Define $X$ and $Y$ q-thrills.}
\begin{proof} %We shall assume $U\subseteq X$ and $V\subseteq Y$  both have size greater than $1/p$ (as the lemma is trivial otherwise).
 
First, assume that $|V| = q |U|$. Let $\F$ be a maximal $X$ $q$-thrill in $G(U,V)$ and let $\F\cap U = \tilde{U}$, i.e., let $\tilde{U}$ denote the set of all those vertices in $U$ which belong to a $q$-fan in $\F$. Similarly, let $\F \cap V = \tilde{V}$ and set $A:=U\setminus \tilde{U}, B:=V\setminus \tilde{V}$. 
Since $\F$ is an $X$ $q$-thrill, $q(u-a) = v-b$ which gives $b = qa$.
Note that we may assume that $a > 1/p$ as otherwise, the bounds on $a$ and $b$ hold trivially since $1/p < d_0/q$ for either assumption on $q$. 
%Again, we shall assume that $a > 1/p$ as the lemma holds trivially otherwise for the same reason.

By the maximality of $\F$, no vertex in $A$ has more than $q-1$ neighbors in $B$, implying $e(A,B) < qa$. Since $a > 1/p$, the aforementioned observation coupled with Theorem \ref{thm:thom} implies
\[
qa > e(A,B) > pab - \sqrt{pnab(1+\varepsilon pa)}
\]
so that
\[pab - \sqrt{pnab(1+\varepsilon pa)} < qa.\]
Plugging $b = qa$ yields
\[
q(pa-1)^2 < pn (1 + \varepsilon pa)
\]
which upon further simplification, yields the following quadratic inequality in $a$:
\begin{equation}\label{eq:quad_a}
qp^2a^2 - (2pq + \varepsilon p^2n) a + q - pn < 0.
\end{equation}
%Note that we have $pn - q  > 0$ for large enough $k,n$ since $q = O_\ep (1)$, and thus,
Since $pn-q>0$ for either assumption on $q$ for $n$ sufficiently large,
\begin{align*}
a &< \frac{2q + \varepsilon pn + \sqrt{(2q+\varepsilon pn)^2 + 4 q( pn - q)}}{2qp}\\
&= \frac{(2q + \varepsilon pn)}{2qp}\left(1+ \sqrt{1+\frac{4q(pn-q)}{(2q + \varepsilon pn)^2}}\right)\\
&<\frac{(2q + \varepsilon pn)}{qp}\left(1 + \frac{2q(pn-q)}{(2q + \varepsilon pn)^2}\right) \text{	(as }\sqrt{1+x}< 1+\frac{x}{2}\text{ for all }x>0)\\
&= \frac{2}{p} + \frac{\varepsilon n}{q} + \frac{2(pn-q)}{p(2q + \varepsilon pn)}<\frac{2}{p} + \frac{\varepsilon n}{q} + \frac{2n}{(2q + \varepsilon pn)}<\frac{2}{p} + \frac{\varepsilon n}{q} + \frac{2}{\varepsilon p} \eqqcolon d.\\
\end{align*}
It now suffices to show that (for either assumption on $q$) $d\leq d_0/q$. Note that $
\frac{2}{p} + \frac{2}{\varepsilon p} < \frac{4k}{\ep \omega(k)}
$.
If $q = O_\ep(1)$, then we have for large enough $k$ that $\omega(k) > 4q/\ep^2$ and therefore, 
%\begin{equation}\label{eq:q1}
$\frac{4k}{\ep \omega(k)} \leq \frac{4n}{\ep \omega(k)}< \frac{\ep n}{q}.
$ %\end{equation}
If $q = \left\lfloor \frac{n}{k}\right\rfloor$, then for large enough $k$, we have that $\omega(k) > 4/\ep^2$ and therefore,
%\begin{equation}\label{eq:q2}
$\frac{4k}{\ep \omega(k)} \leq \ep k\leq \frac{\ep n}{q}.
$%\end{equation}

Now, assume that $u = qv$. This case proceeds analogously to the previous one, with only minor changes at appropriate places. Let $\F$ now be a maximal $Y$ $q$-thrill and let $\tilde{U} = \F\cap U $ and $\F \cap V = \tilde{V}$. Define $A$ and $B$ as in the previous case. Then by the maximality of $\F$, no vertex in $B$ has more than $q-1$ neighbors in $A$, implying $e(A,B) < qb$. Further, we have $a = qb$. By Theorem \ref{thm:thom}, assuming $a > 1/p$ as earlier, we have
\[
qb > pab - \sqrt{pnab(1+\varepsilon pa)}.
\]
Upon plugging in $a = bq$ and working out as before, we obtain the quadratic inequality 
\[
qp^2b^2 - (2pq + \varepsilon p^2 qn)b + (q - pn) < 0
\]
which is identical to (\ref{eq:quad_a}) except with $b$ in place of $a$ and $q\ep$ in place of $\ep$. Thus, it follows that $b < \frac{2}{p} + \varepsilon n + \frac{2}{\varepsilon pq} \leq \frac{2}{p} + \varepsilon n + \frac{2}{\varepsilon p} = d$, therefore $a \leq q d$. This implies the claimed bounds in terms of $d_0$ as before.%by (\ref{eq:q1}) and (\ref{eq:q2}).% and a similar calculation  as in the previous case leads to
%$a < q (\frac{2}{p} + \varepsilon n + \frac{1}{\varepsilon p})$, thus proving $a \leq q d_0$ and $b\leq d_0$.
\end{proof}

A few remarks are in order.
\begin{enumerate}
\item Though we have slightly stronger bounds on $a$ and $b$ in the second case (when $u = qv$), we simply use the stated bounds for the sake of ease of calculations later.
\item When $\ep = 0$ (for instance in the pseudorandom graphs that arise from the point-hyperplane incidences of projective geometries), the calculations above in fact yield 
$a < \frac{1}{p} + \sqrt{\frac{n}{pq}}$ when $v = qu$ and something analogous when $u = qv$. In particular, the sizes of the deleted parts are considerably smaller in this case.
\item If $U\subset X'\subset X, V\subset Y'\subset Y$ then the conclusions of Lemma \ref{lem:main} hold even for the graph $G(X',Y')$  with {\it the same parameters $(p,\ep)$} since the lemma directly applies to the pair $(U,V)$ as a subset of $(X,Y)$. This is vitally of use in the way we apply the Lemma in the proof of Theorem \ref{thm:pseudo} part (b).
\end{enumerate}

%Before we get to the proof of Theorem \ref{thm:pseudo} we make a couple of observations.

\begin{proof}[Proof of Theorem \ref{thm:pseudo} part (a)]
Suppose $n = qk + r$, where $q = \left\lfloor \frac{n}{k} \right\rfloor$ and $r$ is an integer such that $0 \leq r < k$. Choose an arbitrary subset $C_Y\subset Y$ of size $r$ and define $Y_1 = Y\setminus C_Y$. Apply Lemma \ref{lem:main} to the sets $U = X$ and $V = Y_1$ to obtain $A\subset X$ and $B\subset Y_1$ such that $G(X\setminus A, Y\setminus B)$ is spanned by an $X$ $q$-thrill and therefore has NMP (by Lemma \ref{Euclidean:NMP}).  
Define $\X = A$ and $\Y = C_Y \cup B$ so that 
\[\frac{|\X|}{k} \leq \frac{d_0}{qk} \leq 4 \varepsilon = O(\varepsilon)\]
and 
\[\frac{|\Y|}{n} \leq \frac{d_0 + r}{n} < 2\varepsilon + \frac{k}{n} < 3 \sqrt{\varepsilon} = O(\sqrt{\varepsilon}).\]
\end{proof}

\begin{lemma}\label{lemma:smallcase} Suppose $L/\ell$ is representation in reduced form of $n/k$, suppose $L,\ell=O_{\ep}(1)$ and let $d_0 = 2\ep n$. There exist subsets $D_X \subset X,D_Y\subset Y$ with $|D_X|\leq \ell m d_0$ and $|D_Y| \leq L m d_0$, such that $G(X\setminus D_X,Y\setminus D_Y)$ admits a $T_{\ell,L}$-factor. Here, $m$ is the complexity of the Euclidean algorithm  for the parameters $(\ell,L)$  as defined in Section \ref{prelim}. \end{lemma} %is spanned by a vertex-disjoint union of copies of $T_{\ell,L}$. %{\color{red} Could state this better.}
% The assumption that $L,\ell=O_{\ep}(1)$ seems rather restrictive to the general statement of theorem  \ref{thm:pseudo}. But we shall later see that the general case reduces to this special case with only a slight worsening of the bounds. 

\begin{proof} [Proof of Lemma \ref{lemma:smallcase}] Partition both $X$ and $Y$ arbitrarily into ``blocks'', each of size $t = gcd(k,n)$. Let the blocks be denoted by $X_1,\ldots, X_\ell$ and $Y_1,\ldots, Y_L$ respectively. We shall refer to the $X_i$ blocks as \emph{left} blocks and the $Y_j$ blocks as \emph{right} blocks. Let $r_i, q_j$ be the remainders and quotients as defined in Section \ref{prelim}.  We shall now replicate the Euclidean-$(\ell,L)$ process with the vertices being replaced by these blocks, which we shall carry out in $m$ stages, beginning with stage $1$. 

In the rest of the proof of Lemma \ref{lemma:smallcase} we assume that $m$ is even; the $m$ odd case is completely analogous. We also define the sets $\mathcal{X}^{(i)}$ and $\mathcal{Y}^{(i)}$ analogous to the sets $X^{(i)}$ and $Y^{(i)}$ in the definition of the Euclidean tree (see Section \ref{prelim}) as follows.  If $i$ is even,
\[
\mathcal{X}^{(i)} = X_1\sqcup \cdots \sqcup X_{r_{i}}\text{ and } 
\mathcal{Y}^{(i)} = Y_1\sqcup \cdots \sqcup Y_{r_{(i+1)}}
\] 
and if $i$ is odd, then
\[
\mathcal{X}^{(i)} = X_1\sqcup \cdots \sqcup X_{r_{(i+1)}}\text{ and } 
\mathcal{Y}^{(i)} = Y_1\sqcup \cdots \sqcup Y_{r_{i}}
\]
We also assume that $\mathcal{X}^{(0)} = \mathcal{Y}^{(0)} = \emptyset$.

We induct on $m$.  At stage $i$, we apply Lemma \ref{lem:main} to appropriately defined sets $U_{i}$ and $V_{i}$ to obtain sets $A_{i}\subset U_{i}$ and $B_{i}\subset V_{i}$ such that $G(U_{i}\setminus A_{i}, V_{i}\setminus B_{i})$ is spanned by an $X$ $q_i$-thrill or a $Y$ $q_i$-thrill (depending on whether $i$ is even or odd respectively). In fact, it will turn out that $U_{i}$ and $V_{i}$ are large subsets of $\mathcal{X}^{(i)}$ and $\mathcal{Y}^{(i)}\setminus \mathcal{Y}^{(i-1)}$ respectively, when $i$ is even (and something analogous when $i$ is odd).
We denote the set of deleted vertices from $X$ and $Y$ at the \emph{end} of stage $i$ by $D^X_{i}$ and $D^Y_{i}$ respectively, and these are obtained by modifying $A_{i}$ and $B_{i}$ suitably, with the help of $D^X_{i-1}$ and $D^Y_{i-1}$. We then show that $G_i = G(\mathcal{X}^{(i)}\setminus D^X_{i},\mathcal{Y}^{(i)}\setminus D^Y_{i})$ admits a $T_i$-factor, where $T_ i = T_{r_{i},r_{(i+1)}}$ as was defined in Section \ref{prelim}. By controlling the  sizes of $D^X_{i}$ and $D^Y_{i}$ (which we denote by $d^X_{i}$ and $d^Y_{i}$ respectively) the Lemma follows by plugging in $i = m$ because $r_m = \ell$ and $r_{m+1} = L$. 

Let us get to the details now. For starters, we apply Lemma \ref{lem:main} to the ``first'' $r_{1}$ right blocks (recall that $r_{1} = 1$) and the ``first'' $r_{2}$ left blocks. More precisely, we apply Lemma \ref{lem:main} to %$U = \sqcup_{i=1}^{r_{m-2}} X_i$
$U_1 = \mathcal{X}^{(1)} = X_1\sqcup \cdots \sqcup X_{r_{2}}$ and $V_1 = \mathcal{Y}^{(1)} = Y_{r_{1}} = Y_1$ so that $|U_1| = t \cdot r_{2} = t \cdot q_1  r_{1} = q_1 |V_1|$.
We obtain sets $A_1\subset U_1$ and $B_1\subset V_1$ such that $G(U_1\setminus A_1, V_1\setminus B_1)$ is spanned by a $Y$ $q_1$-thrill. This terminates stage $1$ with  $D^X_1:=A_1$ and $D^Y_1:=B_1$; consequently, by Lemma \ref{lem:main} $d^X_1 \leq q_1 d_0$ and $d^Y_1 \leq d_0$. This establishes the following: 

$G_1 = G(\mathcal{X}^{(1)}\setminus D^X_{1}, \mathcal{Y}^{(1)}\setminus D^Y_{1})$ admits a $T_1$-factor, with $d^X_1 \leq q_1 d_0$ and $d^Y_1 \leq d_0$.

  Suppose now that for some $1<i\leq m$, $G_{i-1} = G(\mathcal{X}^{(i-1)}\setminus D^X_{i-1}, \mathcal{Y}^{(i-1)}\setminus D^Y_{i-1})$ admits a $T_{i-1}$-factor,  and 
\begin{enumerate}%[label=(\alph*)]
\item[(1)] if $i$ is even, then  $d^X_{i-1} \leq (i-1)\cdot r_i d_0$ and $d^Y_{i-1} \leq (i-1)\cdot r_{i-1}d_0$. 
\item[(2)]  if $i$ is odd, then  $d^X_{i-1} \leq (i-1)\cdot r_{i-1} d_0$ and $d^Y_{i-1} \leq (i-1)\cdot r_{i}d_0$. 
\end{enumerate}

We shall show that there exist subsets $D^X_i\subset X$ and $D^Y_i\subset Y$ such that $G_i$ admits a $T_i$-factor, and furthermore, 
\begin{enumerate}%[label=(\Alph*)]
\item[(a)] if $i$ is even, then $|D^X_i| = d^X_i \leq i r_i d_0$ and $|D^Y_i| = d^Y_i \leq i r_{i+1} d_0$,
\item[(b)] if $i$ is odd, then $|D^X_i| = d^X_i \leq i r_{i+1} d_0$ and $|D^Y_i| = d^Y_i \leq i r_{i} d_0$,
\end{enumerate}
which would establish the induction step.

Suppose $i$ is even.  %We shall skip the proof of (b)$\Rightarrow$(B) as it is completely analogous. 
Let $S^Y_{i}$ be an arbitrary subset of $Y_{r_{(i-1)} +1} \sqcup\cdots\sqcup Y_{r_{(i+1)}}$ of size $q_{i}\cdot d^X_{i-1}$. Define
\[
U_{i} := \mathcal{X}^{(i)} \setminus D^X_{i-1} \text{ and }
V_{i} := (\mathcal{Y}^{(i)}\setminus \mathcal{Y}^{(i-1)})\setminus S^Y_{i} = (Y_{r_{(i-1)} +1} \sqcup\cdots\sqcup Y_{r_{(i+1)}}) \setminus S^Y_{i}
\] 
 Since  $r_{i+1} - r_{i-1} = q_ir_i$ we have $|V_{i}| = t(r_{i+1} - r_{i-1}) - q_id^X_{i-1} =  q_{i} |U_{i}|$, so  by  Lemma \ref{lem:main}, we obtain sets $A_i \subset U_i$ and $B_i\subset V_i$  with $|A_i| \leq d_0/q_i$ and $|B_i| \leq d_0$ such that $G(U_{i}\setminus A_{i}, V_{i}\setminus B_{i})$ is spanned by an $X$ $q_i$-thrill. 

By assumption, $G_{i-1}$ admits a $T_{i-1}$-factor i.e., $G_{i-1}$ is spanned by vertex-disjoint copies of $T_{i-1}$.  Define $\mathrm{CORRUPT}^X_i$ to be the set of all those vertices in $\mathcal{X}^{(i-1)}\setminus D^X_{i-1}$ which belong to one of the above copies of $T_{i-1}$ that also contains at least one vertex from $A_i$. Obviously, $A_i \subseteq \mathrm{CORRUPT}^X_i$. Similarly, we define $\mathrm{CORRUPT}^Y_i$ as the set of vertices in $\mathcal{Y}^{(i-1)}\setminus D^Y_{i-1}$ which belong to a copy of $T_{i-1}$ that contains at least one vertex from $A_i$. We refer to such copies of $T_{i-1}$ in $G_{i-1}$ (that contain at least one vertex from $A_i$) as \emph{corrupt} copies. Define
\[\mathrm{CORRUPT}_i := \mathrm{CORRUPT}^X_i \sqcup \mathrm{CORRUPT}^Y_i\] 
as the set of those vertices of $G_{i-1}$ that get ``corrupted'' due to the introduction of further deletions during stage $i$ (i.e. the set $A_i$). In other words, $\mathrm{CORRUPT}_i$ is the set of vertices touched by the corrupt copies. See Figure \ref{fig:ind_step} for an illustration of the induction step.

\begin{figure}[ht]
\centering
\includegraphics[scale=0.8]{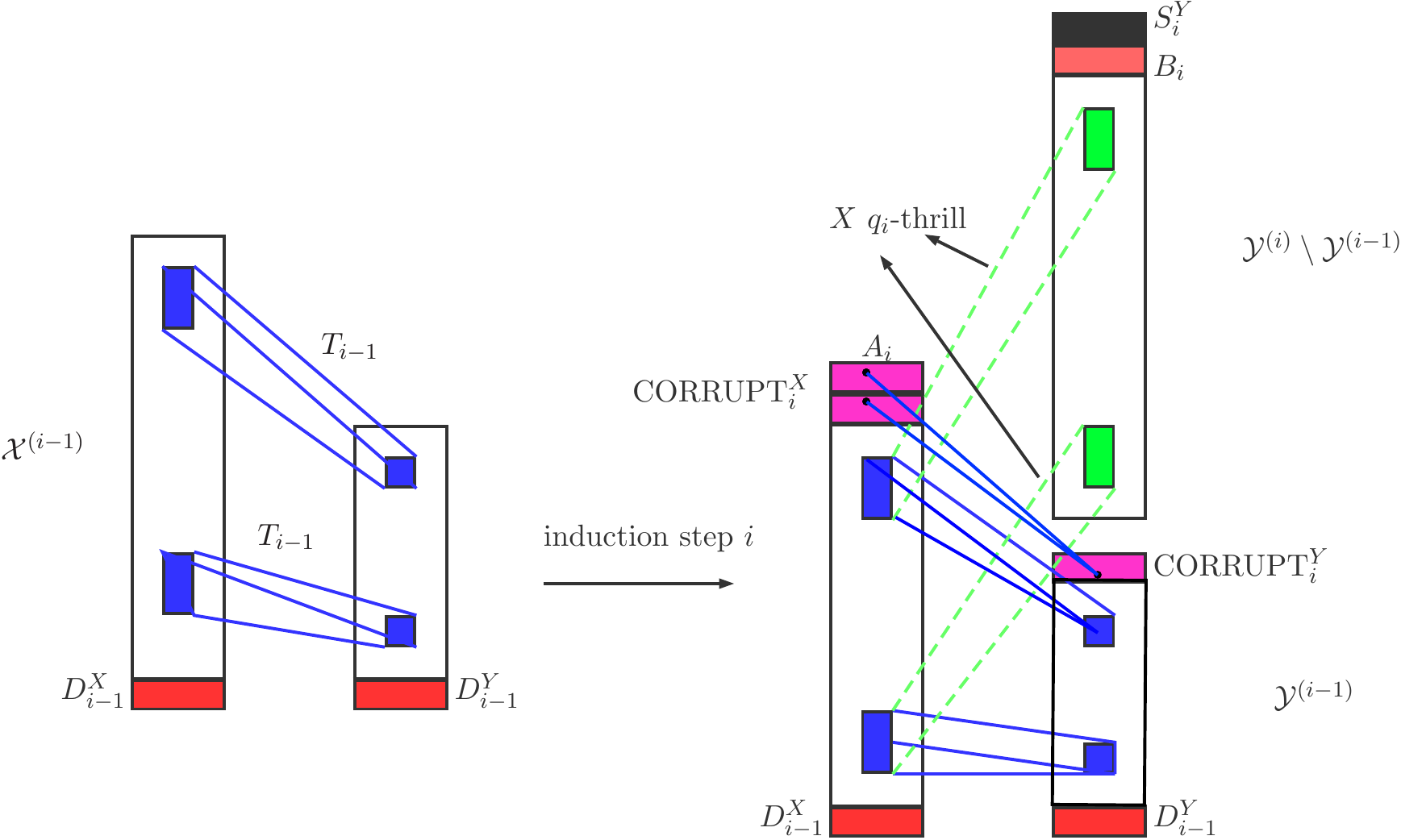}
\caption{An illustration of the induction step in the proof of Theorem \ref{thm:pseudo}. The picture on the left depicts the copies of $T_{i-1}$ that span $G_{i-1}$ and are colored blue. The picture on the right depicts what happens to each of these copies in the induction step: those which have a vertex in $A_i$ (the topmost box in $\mathcal{X}^{(i)}$) ``corrupt'' \emph{all} the vertices that they contain (colored pink) and those which do \emph{not} have a vertex in $A_i$ ``evolve'' to $T_i$ via an $X$ $q_i$-thrill into $\mathcal{Y}^{(i)}\setminus \mathcal{Y}^{(i-1)}$, shown in green.}\label{fig:ind_step}
\end{figure}

%Next, we shall describe the sets $D^X_i$ and $D_i^Y$ using $A_i, B_i$ and the above and estimate their sizes. Formally, we 
Define 
\[
D^X_i := D^X_{i-1} \sqcup \mathrm{CORRUPT}^X_i \text{ and } D^Y_i \coloneqq D^Y_{i-1} \sqcup 
S^Y_i \sqcup B_i \sqcup \mathrm{CORRUPT}^Y_i
\]
and set $d^X_i\coloneqq|D^X_i|, d^Y_i\coloneqq|D^Y_i|$. 
Note that every corrupt copy of $T_{i-1}$ in $G_{i-1}$ has $r_i$ vertices in $X$ and $r_{i-1}$ vertices in $Y$. Therefore, we have the bounds 
\[
|\mathrm{CORRUPT}^X_i| \leq r_i |A_i| \leq \frac{r_i}{q_i}d_0
\text{ and }
|\mathrm{CORRUPT}^Y_i| \leq r_{i-1} |A_i| \leq \frac{r_{i-1}}{q_i} d_0
\]
Putting things together, we obtain the recurrences 
\[
d^X_i \leq d^X_{i-1} + \frac{r_i}{q_i}d_0
\text{ and }
d^Y_i \leq d^Y_{i-1} + q_i d^X_{i-1} + d_0 + \frac{r_{i-1}}{q_i} d_0
\]
By the induction hypothesis we have $d^X_{i-1} \leq (i-1)\cdot r_i d_0$ and $d^Y_{i-1} \leq (i-1)\cdot r_{i-1}d_0$. Therefore, 
\[
d^X_i \leq d_0\left((i-1)\cdot r_i  + \frac{r_i}{q_i}\right) \leq i\cdot r_i d_0
\]
and 
\[
d^Y_i \leq d_0\left( (i-1)\cdot r_{i-1} + (i-1)\cdot q_i r_i + 1 + \frac{r_{i-1}}{q_i}\right) \leq i \cdot r_{i+1} d_0
\]
where in the final step, we use $r_{i+1} - r_{i-1} = q_ir_i$ and the fact that $1 + r_{i-1} \leq r_i < r_{i+1}$. 

We now prove that $G_i$ admits a $T_i$-factor. 
Recall from the preliminaries that if $T_{i-1} = T_{r_i,r_{(i-1)}}$ is the Euclidean tree with \emph{left} vertices $x_1,\ldots,x_{r_i}$ and \emph{right} vertices $y_1,\ldots, y_{r_{(i-1)}}$, then $T_i = T_{r_i,r_{(i+1)}}$ is constructed on \emph{left} vertices $x_1,\ldots,x_{r_i}$ and \emph{right} vertices $y_1,\ldots, y_{r_{(i+1)}}$, by adding to $T_{i-1}$ an $X$ $q_i$-thrill of size $r_i$ between $x_1,\ldots,x_{r_i}$ and $y_{r_{(i-1)}+1},\ldots, y_{r_{(i+1)}}$. 
By Lemma \ref{lem:main}, $G(\mathcal{X}^{(i)}\setminus D^X_{i},(\mathcal{Y}^{(i)}\setminus \mathcal{Y}^{(i-1)} )\setminus D^Y_{i})$ is spanned by an $X$ $q_i$-thrill. 
This, along with the copies of $T_{i-1}$ that span $G_{i-1}$, gives us the desired $T_i$-factoring of $G_i$. %$G(\mathcal{X}^{(i)}\setminus D^X_{i},\mathcal{Y}^{(i)}\setminus D^Y_{i})$.

The proof of the inductive step when $i$ is odd i.e., $(2)\Rightarrow(b)$ is completely analogous ($X$ swapped with $Y$ everywhere). The only small difference that arises is in the recurrences for $d^X_i$ and $d^Y_i$ because of the slightly different bounds for $|A_i|$ and $|B_i|$ given by Lemma \ref{lem:main} in this case. In particular, by following the same line of argument as in the proof of $(1)\Rightarrow(a)$, we obtain, in this case 
\[
d^X_i \leq d^X_{i-1} + q_id^Y_{i-1} + q_i d_0 +  r_{i-1} d_0 
\text{ and }
d^Y_i \leq d^Y_{i-1} + r_i d_0.
\]
But then, by using the trivial bound $q_i + r_{i-1} \leq r_{i+1}$, we obtain the desired estimates $d^X_i \leq i\cdot r_{i+1} d_0$ and $d^Y_i \leq i\cdot r_i d_0$.

Thus, we have shown that there exist subsets $D_X = D^X_m \subset X,$ $D_Y = D^Y_m\subset Y$ such that $G(X\setminus D_X,Y\setminus D_Y)$ admits a $T_{\ell,L}$-factor and consequently has NMP. Furthermore we have $$|D_X|\leq \ell m d_0\textrm{\ and\ } |D_Y| \leq L m d_0.$$ \end{proof}
%These estimates are clearly going to be more helpful when we have a \emph{small} representation in reduced from of the fraction $n/k$ i.e., when $n$ and $k$ have a large gcd.  
%Recall that we had used the assumption that $L \leq \frac{1}{\sqrt[4]{\varepsilon^3}}$ crucially in the proof of Lemma \ref{lem:main} (not the specific bound itself, but simply that $L$ is bounded by a growing function of $1/\varepsilon$).

%Next, we show how we can always achieve so by simply considering a suitably defined \emph{large} subgraph of the original Thomason pseudorandom graph given to us. More formally, 

We are now in a position to prove Theorem \ref{thm:pseudo} part (b). %converts any Thomason pseudorandom graph to one with part sizes whose ratio has a ``small'' representation in lower terms, as was assumed in Lemma \ref{lemma:smallcase}.
\begin{proof}[Proof of Theorem \ref{thm:pseudo} part (b)] Suppose $G$ is a Thomason pseudorandom bipartite graph with parameters $(p,\varepsilon)$ and with vertex classes $X$ and $Y$ of sizes $k$ and $n$  respectively with $\frac{n}{k} \leq \frac{1}{\sqrt{\varepsilon}}$.

Set $\alpha \coloneqq \sqrt[4]{\varepsilon^3}$ and $\eta \coloneqq \sqrt[4]{\varepsilon}$ and consider the interval $[n(1-\alpha),n]$. Since its length is $\alpha n$, there is an integer $N\in I$ such that $N$ is a multiple of $\lfloor \alpha n \rfloor$. Also, since $\eta k \geq \alpha n$, there is an integer $K$ in the interval $J = [k(1-2\eta), k(1-\eta)]$ such that $K$ is a multiple of $\lfloor \alpha n \rfloor$. With $K$ and $N$ as defined above (note that $K\leq N$), simply pick a subset $C_X\subset X$ of size $k-K$ and $C_Y \subset Y$ of $n-N$ arbitrarily and define a new graph $G' = G(X\setminus C_X, Y\setminus C_Y)$. Observe that if $L/\ell$ is the representation in reduced from of  $N/K$, then $L\leq \frac{1}{\sqrt[4]{\varepsilon^3}} = O_\ep(1)$. Applying Lemma \ref{lemma:smallcase} to $G'$ (see Remark $3$ after Lemma \ref{lem:main}), 
%we repeat the entire argument above to
we obtain subsets $D_X \subset X\setminus C_X$ and $D_Y \subset Y\setminus C_Y$ such that $G(X\setminus \X, Y\setminus \Y)$ has NMP, where $\X = C_X\cup D_X$ and $\Y = C_Y\cup D_Y$. By Fact \ref{fact:euclid} and the trivial bounds $K\leq n$ and $N\leq n$, we have
\[
\frac{|\Y|}{n} \leq \alpha + \frac{L}{n}\cdot m d_0 \leq \sqrt[4]{\varepsilon^3} + 5 \sqrt[4]{\varepsilon}\log L \leq 6\sqrt[4]{\varepsilon}\log \left(\frac{1}{\varepsilon}\right)% = O(\varepsilon^{1/5})
\]
and similarly, 
\[
\frac{|\X|}{k} \leq 2\eta + \frac{\ell}{k} \cdot m d_0 \leq 7 \sqrt[4]{\varepsilon}\log \left(\frac{1}{\varepsilon}\right)% = O(\varepsilon^{1/5}).
\] and that completes the proof. 
\end{proof}

\section{Concluding Remarks}\label{sec:conc}
\begin{itemize}
\item The main engine in the proof of Theorem  \ref{thm:pseudo} comes from Lemma \ref{lem:main} which is the place the pseudorandomness is used in an explicit form. The rest of the proof of the theorem including the inductive argument uses this in a black-box manner. Hence, if we had an equivalent statement to Lemma 4.1 for {\it other} models of pseudorandomness - call it Lemma 4.1* (say), then the rest of the proof of Theorem 1.4 can run through with the error estimates being dictated by Lemma 4.1* instead. The content of Lemma \ref{lem:main} uses the notion of Thomason pseudorandomess explicitly {\it only} when we evoke Theorem \ref{thm:thom} which is basically a statement that estimates how much the difference between $e(A,B)$ and the expected number of edges, if the graph were random, viz., $p|A||B|$ can be. For $(n,d,\lambda)$ graphs, the analogue of this theorem is the expander-mixing lemma which provides precisely such an estimate. 

%\item As the hypothesis of Lemma \ref{lem:main} only requires the applicability of Theorem \ref{thm:thom} which is essentially an expander-mixing type lemma, {\it any} model of pseudorandom bipartite graphs with a corresponding expander-mixing theorem can be invoked to get a result similar to that of Lemma \ref{lem:main}. 
We illustrate this by returning to problem \ref{prob:sumcayley} that was stated in the introduction. For $\ep>0$, and $q$ a sufficiently large prime power, let $H$ be a multiplicative subgroup of $\bF_q^*$ of order at least $q^{1/2+\ep}$. Consider the Sum-Cayley graph $\Gamma_q(H)$ whose vertex set is $\bF_q$ and vertices $x, y$ are adjacent if and only if $x+y\in H$.  A result of Alon and Bourgain (see \cite{AB})  states that that $\Gamma_q(H)$ is a $(q,|H|,q^{1/2})$ graph, i.e., it is a regular graph on $q$ vertices, with degree $|H|$, and every non-trivial eigenvalue of $\Gamma_q(H)$ is at most $q^{1/2}$.  If $G$ is  the bipartite graph described in the introduction following the description of problem 2, then it is not difficult to show that for any $A\subset X, B\subset Y$ we have $|e(A,B) - \frac{|A||B||H|}{q}|<\sqrt{q|A||B|}$ by using the expander-mixing lemma.  Then, via the argument in the proof of  Lemma \ref{lem:main} we have: If $X, Y\subset\bF_q$ with $|Y|=10|X|$, $|X|\ge q/100$, and let $H$ is a subgroup of $\bF_q^*$ of size at least $q^{1/2+\ep}$, then there exists $A\subset X, B\subset Y$ with $|A|\le O(q^{1-\ep})$, and $|B|=10|A|$ such that $G(X\setminus A, Y\setminus B)$ has NMP. Consequently, every element of $Y\setminus B$ can be labeled by some element of $X\setminus A$ such that  each label appears $10$ times, and further, for each $y\in Y$ labeled $x$, the sum $x+y\in H$. This answers in the affirmative, the approximate version of problem 2. One could pose more general questions of the same kind, but without the additional constraint that $|Y|$ is a multiple of $|X|$. For instance, suppose $X,Y\subset\bF_q$ and $|Y|=\frac{3}{2}|X|$ (say), with $|X|\ge \Omega(q)$, and let $H$ be a subgroup of $\bF_q^*$ of size at least $q^{1/2+\ep}$. Then one can similarly show that there exist subsets $\X\subset X, \Y\subset Y$ with $|\X|\le f(\ep)|X|, |\Y|\le g(\ep)|Y|$ such that if $X', Y'$ are the remaining sets, then one may form a star-array $\A$ of dimension $|X'|\times |Y'|$ whose rows and columns are labeled by the elements of $X', Y'$ respectively with the property that if the $(x,y)^{th}$ element of $\A$ is a star, then $x+y\in H$. Furthermore, each row of $\A$ has precisely $3$ stars, and each column has precisely $2$ stars. 
\item For a bipartite graph $G(X,Y)$ with $|X|=|Y|$ that admits a perfect matching, the Max-Min Greedy Matching problem that was introduced in \cite{eff} goes as follows. Given permutations $\sigma, \pi$ of the vertices of $X$ and  $Y$ respectively, the vertices of $X$ are processed according to  $\sigma$, and each $x\in X$ is matched to its earliest available neighbor in $Y$ according to  $\pi$. If $M_G[\sigma,\pi]$ denote the size of the resulting greedy matching,  determine $\rho[G]\coloneqq\frac{\max_{\pi}\min_{\sigma} |M_G[\sigma,\pi]|}{|X|}$. This problem admits a natural generalization. Suppose $G(X,Y)$ is a bipartite graph, with $|X|=k, |Y|=n$, with $k\le n$, and suppose $r=\lfloor n/k\rfloor$. As before, let $\sigma, \pi$ be permutations of the vertices of $X$ and $Y$ respectively. We process the vertices of $X$ according to $\sigma$ and for each $x\in X$, we choose its {\it first $r$ neighbors} in $Y$ that have not been already chosen by some previous vertex of $X$ according to $\pi$.  Let $m^{(r)}_G[\sigma,\pi]$ denote the number of vertices of $X$ for which one can choose $r$ such neighbors. Then determine $\rho_r[G]\coloneqq\frac{\max_{\pi}\min_{\sigma} m^{(r)}_G[\sigma,\pi]}{|X|}$. Our proof of Lemma \ref{lem:main} can easily be adapted to establish the following: Suppose $\ep>0$, and let $\omega$ be a function such that $\omega(k)\to\infty$ as $k\to\infty$. Then there exists $k_0=k_0(\ep)$ such that whenever $n\ge k> k_0$ and $G(X,Y)$ is a  $(p,\ep)$-Thomason pseudorandom bipartite graph with $|X|=k, |Y|=n$, and  $p\ge\frac{\omega(k)}{k}$, then  $\rho_r[G]\ge 1-O(\ep)$. %We omit the details.
\item Our proof of Theorem \ref{thresh_NMP} on closer examination reveals that $\bG(k,n,p)$ does not have NMP {\it whp} for $p=\frac{\log n-\omega(n)}{k}$ for any arbitrary function $\omega$ that goes to infinity.  However, to prove the existence of NMP with high probability,  our proof cannot extend beyond $p=\frac{\log n+O(\sqrt{\log n})}{n}$. While it is possible to improve (using our methods) our result to prove that $\bG(k,n,p)$ has NMP {\it whp} for $p=\frac{\log n + f(n)}{k}$ for some $f=o(\log n)$, the question of whether there is a sharp threshold for NMP of the form $p=\frac{\log n+\omega(n)}{k}$ remains open.   
\item As remarked in the Introduction, our proof of Theorem \ref{thm:pseudo} shows that $f(x)=g(x)=O(x^{1/4}\log(1/x))$ works uniformly for all pairs $(k,n)$. Is it possible to improve this to $f(x)=g(x)=O(x)$ uniformly over all $(k,n)$? 
\item We make a  final remark pertaining to a remark following the statement of Theorem \ref{thm:pseudo} in the Introduction. As we noted, the definition of Thomason pseudorandomness does not preclude the existence of isolated vertices unless a more symmetric definition of pseudorandomness is adopted. In that case, it would be interesting to see if one can arrive at a stronger conclusion than the statement of Theorem \ref{thm:pseudo}. 
\end{itemize}

% \newpage
%\section*{Appendix: Robustness of Thomason pseudorandomness}
\begin{appendices}
\newpage
\section*{Appendix: Robustness of Thomason pseudorandomness}

\paragraph*{Lemma.}
{\it
Let $0<\varepsilon <\frac{1}{2}$, and $k\le n$ be positive integers. Suppose $G(X,Y)$ is a Thomason pseudorandom bipartite graph with parameters $(p_0,\varepsilon_0)$ with $|X|=k, |Y|=n$, and suppose $p_0\geq \frac{1}{\sqrt{k}}$. Then, for a given integer $D$ satisfying $\frac{\alpha}{2} n\leq D\leq \alpha n$ for $\alpha = \varepsilon^3$, there exist subsets $C_X\subseteq X$ and $C_Y\subseteq Y$ such that 
\begin{itemize}
\item $|C_Y| = D$ and $|C_X| \leq \eta k$, where $\eta = 2\exp(-\frac{C}{\ep})$ for some fixed constant $C$,
\item the subgraph induced by the sets $X\setminus C_X$ and $Y\setminus C_Y$ is Thomason pseudorandom with parameters $(p_1,\varepsilon_1)$ where $p_1 = p_0 (1-\varepsilon)$ and $\varepsilon_1\le 5 (\ep_0 + 3\ep)$. %\frac{\ep_0+2\ep+\alpha}{(1-\varepsilon)^2(1-\alpha)}$
\end{itemize}
}
\begin{proof} Let $\eta=2\exp(-\frac{C}{\ep})$ where $C$ shall be specified later. Let $T\subseteq Y$  be a uniformly random subset of $Y$  of size $D$. Then by the tail bound of the hypergeometric distribution (see \cite{skala}) we have, for every $t\geq 0$, 
\begin{equation}\label{eqn:hgtail}
\Pr\left[\left||N(u)\cap T| -\frac{d(u)}{n}D\right| \geq tD\right] \leq 2e^{-2t^2D}
\end{equation}
 for every vertex $u\in X$. Now, fix $t = \varepsilon p_0 (\frac{n}{D}-1)$. Call a vertex $u\in X$ \emph{bad} with respect to $T$ if 
\[
|N(u)\cap T|\geq \left(\frac{d(u)}{n} + t\right)D.
\]
Then by equation~\ref{eqn:hgtail}, the expected number of bad vertices is at most $2k e^{-2t^2D}$. Fix a set $C_Y\subseteq Y$ of size $D$ for which the set of bad vertices (which we shall call $C_X$) has size at most $2k e^{-2t^2D}$. 

Now, for a vertex $x\in X$, let $N'(x) = N(x)\cap (Y\setminus C_Y)$. Then for $x\in X\setminus C_X$, as $x$ is not a bad vertex,
\[
|N'(x)| = |N(x)| - |N(x)\cap C_Y|  \geq p_0n\left(1-\frac{D}{n}\right) - Dt = p_0(1-\varepsilon)(n-D) = p_1\cdot |Y\setminus C_Y|
\]
where the inequality follows from the hypothesis (see Definition~\ref{def:pseudo}) that $G$ is Thomason pseudorandom.   Also note that for any distinct vertices $u,v\in X\setminus C_X$, 
\[
|N'(u)\cap N'(v)| \leq (1+\varepsilon_1)p_1^2\cdot|Y\setminus C_Y|\cdot
\]
which follows since 
\[|N'(u)\cap N'(v)|\le |N(u)\cap N(v)|\le (1+\varepsilon_0)p_0^2 n  \leq \left(\frac{1+\ep_0}{(1-\varepsilon)^2(1-\alpha)}\right)p_1^2 (n-D) \]
where the  last inequality follows from the fact that  $n-D \geq n(1-\alpha)$. The required codegree bound then follows from the given condition on $\ep_1$.

It remains only to check is that $2e^{-2t^2D}\leq\eta$. To see this, observe that $t = \varepsilon p_0 (\frac{n}{D}-1) \geq \varepsilon p_0 (\frac{1}{\alpha}-1)$ and also note that $\varepsilon<\frac{1}{2}\Rightarrow 1-\alpha > \frac{7}{8}$. Thus,
\[
2t^2D \geq 2\varepsilon^2 p_0^2 \left(\frac{1}{\alpha}-1\right)^2 \left(\frac{\alpha n}{2}\right)\geq \left(\frac{\varepsilon^2}{\alpha}\right)(1-\alpha)^2(p_0^2k)\geq \frac{49}{64\varepsilon}= \log \left(\frac{1}{\eta}\right)
\]
where we may take the constant $C=\frac{49}{64}=0.765625$ in the definition of $\eta$.
\end{proof}
One interesting consequence of the proof of the lemma is that if we seek $\eta=\textrm{poly}(\ep)$ then one has a randomized algorithm to choose a set $T\subset Y$ and a related $BAD(T)\subset X$ with $|T|=D, |BAD(T)|\le\eta k$ such that deleting these sets from $Y, X$ respectively results in another Thomason pseudorandom graph with only slightly worse parameters.

It is known (see \cite{thom2}) that bipartite graphs arising from the point-hyperplane incidence structure of a projective geometry of dimension $d$ over a finite field $\bF_q$  is Thomason pseudorandom with parameters $p=n^{-1/2}(1+o(1))$ and $\ep=0$. More generally, one can take the point-block incidence structure arising from a symmetric block design as the ``seed'' Thomason pseudorandom graph which upon the application of the lemma above gives us several other examples of Thomason pseudorandom graphs with parameters that are relevant in Theorem \ref{thm:pseudo}.
\end{appendices}
\end{document}